\documentclass[11pt]{article}
\usepackage{color}
\usepackage{amssymb,latexsym}

\newtheorem{lem}{Lemma}[section]
\newtheorem{theo}{Theorem}[section]
\newtheorem{pro}{Proposition}[section]
\newtheorem{cor}{Corollary}[section]
\newtheorem{con}{Conjecture}[section]

\newcommand {\red} {\textcolor{red}}
\newcommand {\green} {\textcolor{green}}
\newcommand {\blue} {\textcolor{blue}}

\renewcommand{\theenumi}{\rm (\roman{enumi})}

\newcommand{\proof}
{{\noindent {\em Proof}.\quad}\setcounter{countclaim}{0}
\setcounter{countcase}{0}}
\newcommand{\proofend}{{\hfill$\Box$}}
\newcommand{\proofends}{\mbox {\hfill    } \Box}
\newcommand{\gap}{\vspace{0.5 cm}}
\newcommand{\sgap}{\vspace{0.3 cm}}

\newcommand{\rem}{\noindent {\bf Remark.}}
\newcommand{\rems}{\noindent {\bf Remarks.}}
\newcommand{\note}{\noindent {\bf Note.}}
\newcommand{\notes}{\noindent {\bf Notes.}}

\newcounter{countfig}
\def\infig{\addtocounter{countfig}{1}{Figure \thecountfig}}
\def\itfig{{\rm {\it Figure} \thecountfig}}
\def\norfig{{\rm Figure \thecountfig}}

\newcounter{countclaim}
\def\inclaim{\addtocounter{countclaim}{1}
{\noindent {\bf Claim \thecountclaim}: }}

\newcounter{countcase}
\def\incase{\addtocounter{countcase}{1}
{\noindent {\bf Case \thecountcase}: }}

\newcommand{\beeq}{\begin{equation}}
\newcommand{\eneq}{\end{equation}}

\newcommand{\beeqn}{\begin{eqnarray*}}
\newcommand{\eneqn}{\end{eqnarray*}}

\def \sete{{\cal E}}
\def \setf{{\cal F}}
\def \setsf{{\cal SF}}
\def \setp{{\cal P}}
\def \sets{{\cal S}}
\def\setz{{\cal Z}}
\def \setr{{\cal R}}
\def \sett{{\cal T}}
\def \setst{{\cal ST}}
\def \setec{{\cal C}}

\def\relabel {\label}  %final version

\begin{document}

\newcommand{\resection}[1]
{\section{#1}\setcounter{equation}{0}}

\renewcommand{\theequation}{\thesection.\arabic{equation}}

\renewcommand{\theenumi}{\rm (\roman{enumi})}
\renewcommand{\labelenumi}{\rm(\roman{enumi})}

\baselineskip 0.6 cm

\title {Expression for the Number of Spanning Trees of 
Line Graphs of Arbitrary Connected Graphs\thanks{This paper was partially supported by 
NSFC (No. 11271307, 11171134 and 11571139)
and NIE AcRf (RI 2/12 DFM) of Singapore.
}
}

\author
{Fengming Dong\thanks{Corresponding author.
Email: fengming.dong@nie.edu.sg}\\
\small Mathematics and Mathematics Education\\
\small National Institute of Education\\
\small Nanyang Technological University, Singapore 637616
\sgap \\
Weigen Yan\\
\small School of Sciences, 
Jimei University, Xiamen 361021, China}

\date{}

\maketitle

\begin{abstract}
For any graph $G$, 
let $t(G)$ be the number of spanning trees of $G$, 
$L(G)$ be the line graph of $G$ and 
for any non-negative integer $r$, 
$S_r(G)$ be the graph obtained from $G$ by
replacing each edge $e$ by a path of length $r+1$ connecting
the two ends of $e$. 
In this paper we 
obtain an expression for $t(L(S_r(G)))$ 
in terms of spanning trees of $G$ by a combinatorial approach.
This result generalizes some known results on the relation between $t(L(S_r(G)))$ and $t(G)$
and gives an explicit expression  
$t(L(S_r(G)))=k^{m+s-n-1}(rk+2)^{m-n+1}t(G)$
if $G$ is of order $n+s$ and size $m+s$ in which
$s$ vertices are of degree $1$ and the others are
of degree $k$. 
Thus we prove a conjecture on 
$t(L(S_1(G)))$ for such a  graph $G$.
\end{abstract}


\resection{Introduction}

The graphs considered in this article 
have no loops but may have parallel edges.
For any graph $G$,
let $V(G)$ and $E(G)$ be the vertex set and edge set of $G$ respectively,  
let $S(G)$ be the graph obtained from $G$ by inserting a new
vertex to each  edge in $G$, 
$L(G)$ be the line graph of $G$,
$\sett(G)$ be the set of spanning trees of $G$ and 
$t(G)=|\sett(G)|$.
Note that for any parallel edges $e$ and $e'$ in $G$,
$e$ and $e'$ are two vertices in $L(G)$ 
joined by two parallel edges. 
For any disjoint subsets $V_1, V_2$ of $V(G)$, 
let $E_G(V_1, V_2)$ (or simply $E(V_1, V_2)$) denote the set 
of those edges in $E(G)$ which have ends in $V_1$ and $V_2$ 
respectively, and let 
$E_G(V_1, V(G)-V_1)$ be simply denoted by $E_G(V_1)$.
For any  $u\in V(G)$, 
let $E_G(u)$ (or simply $E(u)$) denote the set $E_G(\{u\})$.
So the degree of $u$ in $G$, denoted by $d_G(u)$ (or simply $d(u)$), is equal to $|E(u)|$.
For any subset $U$ of $V(G)$, 
let $G[U]$ denote the subgraph of $G$ induced by $U$
and let $G-U$ denote the subgraph of $G$ induced by $V(G)-U$.
For any $E'\subseteq E(G)$, 
let $G[E']$ be the spanning subgraph of $G$ with edge set $E'$,
$G-E'$ be the graph $G[E(G)-E']$
and $G/E'$ be the graph obtained from 
$G$ by contracting all edges of $E'$. 

Our paper concerns the relation between $t(G)$ and $t(L(G))$ 
or  $t(L(S(G)))$.
Such a relation was 
first found by Vahovskii~\cite{vah}, then by Kelmans~\cite{kel} 
and was rediscovered by Cvetkovi\'{c}, Doob and Sachs~\cite{cve}
for regular graphs.
They showed that 
if $G$ is a $k$-regular graph of order $n$ and size $m$, then
\begin{equation}\relabel{eq1-1}
t(L(G))=k^{m-n-1}2^{m-n+1}t(G).
\end{equation}
The first result on the relation between $t(G)$ and $t(L(S(G)))$
was found by Zhang, Chen and Chen \cite{zha}.
They proved that if $G$ is $k$-regular, then 
\begin{equation}\relabel{eq1-2}
t(L(S(G)))=k^{m-n-1}(k+2)^{m-n+1}t(G).
\end{equation}
Yan~\cite{yan1} recently generalized the result
of (\ref{eq1-1}). 
He proved that if $G$ is a graph of order $n+s$ and size $m+s$ 
in which $s$ vertices are of degree $1$ and all others 
are of degree $k$, where $k\ge 2$, then 
\begin{equation}\relabel{eq1-3}
t(L(G))=k^{m+s-n-1}2^{m-n+1}t(G).
\end{equation}
Yan~\cite{yan1} also proposed a conjecture to generalize the 
result of (\ref{eq1-2}). 

\begin{con}[\cite{yan1}]\relabel{con-yan}
Let $G$ be a connected graph of order $n+s$ and size $m+s$ 
in which $s$ vertices are of degree $1$ and all others 
are of degree $k$.
Then
$$
t(L(S(G)))=k^{m+s-n-1}(k+2)^{m-n+1}t(G).
$$
\end{con}

If $G$ is a digraph, 
the relation between $t(G)$ and $t(L(G))$ 
was first obtained by Knuth \cite{knu} 
by an  application of the Matrix-Tree Theorem
and a bijective proof of the result 
was found by Bidkhori and Kishore \cite{bid}. 
Note that expressions  
(\ref{eq1-1}),  (\ref{eq1-2}) and (\ref{eq1-3})
were also obtained by the respective authors 
mentioned above by an application of the Matrix-Tree Theorem.
To our knowledge, these results still do not have  
any combinatorial proofs.  
Some related results can be seen in 
\cite{ber, che, lev, per, zha2}.

For an arbitrary connected graph $G$
and any non-negative integer $r$, let 
$S_r(G)$ denote the graph obtained from $G$ by
replacing each edge $e$ of $G$ by a path of length $r+1$ connecting the two ends of $e$.
Thus $S_0(G)$ is $G$ itself and $S_1(G)$ is the graph $S(G)$.
Our main purpose in this paper is 
to use a combinatorial method 
to find an expression for $t(L(S_r(G)))$ %which is 
given in Theorem~\ref{gen-Sr}.

\begin{theo}\relabel{gen-Sr}
For any connected graph $G$ and any integer $r\ge 0$, 
\begin{equation}\relabel{eq-r-ano}
t(L(S_r(G)))=\prod_{v\in V(G)}d(v)^{d(v)-2}
\sum_{E'\subseteq E(G)}t(G[E'])
r^{|E'|-|V(G)|+1}
\prod_{e\in E(G)-E'}(d(u_{e})^{-1}+d(v_{e})^{-1}),
\end{equation}
where $d(v)=d_G(v)$ and $u_e$ and $v_e$ are the two ends of $e$.
\end{theo}

As $S_0(G)$ is $G$ itself, 
the following expression for $t(L(G))$ is a
special case of Theorem~\ref{gen-Sr}:
%For $r=0$, (\ref{eq-r-ano}) is 
\begin{equation}\relabel{eq-r0-ano}
t(L(G))=\prod_{v\in V(G)}d(v)^{d(v)-2}
\sum_{T\subseteq \sett(G)}
\prod_{e\in E(G)-E(T)}(d(u_{e})^{-1}+d(v_{e})^{-1}).
\end{equation}

The proof of Theorem~\ref{gen-Sr} 
will be completed in Sections 3 and 4. 
%for the case $r=0$ (i.e.,the result (\ref{eq-r0-ano})) and in Section 4 for the case $r\ge 1$ respectively.
In Section 3, we will show that 
the case $r=0$ of Theorem~\ref{gen-Sr} 
(i.e., the result (\ref{eq-r0-ano}))
is a special case 
of another result
(i.e., Theorem~\ref{gen-result}), and 
in Section 4, we will prove 
the case $r\ge 1$ of Theorem~\ref{gen-Sr} 
by applying this theorem 
for the case $r=0$ (i.e., (\ref{eq-r0-ano})).
%\ref{eq-r0-ano}).
To establish Theorem~\ref{gen-result},
we need to apply a result 
in Section 2 (i.e., Proposition~\ref{no-trees-original}),
which determines the number of  spanning trees in 
a graph $G$ with a clique $V_0$ 
such that $F=G-E(G[V_0))$ is a forest and 
every vertex in $V_0$ is incident with 
at most one edge in $F$.
Finally, in Section 5, we will apply Theorem~\ref{gen-Sr} to 
show that for any graph $G$ mentioned in Conjecture~\ref{con-yan}
and any  integer $r\ge 0$, we have 
\begin{equation}\relabel{gen-con-yan}
t(L(S_r(G)))=k^{m+s-n-1}(rk+2)^{m-n+1}t(G).
\end{equation}
Thus (\ref{eq1-3}) follows 
and Conjecture~\ref{con-yan} is proved.

Note that in the proof of Theorem~\ref{gen-Sr}, 
we will express $t(L(S_r(G)))$ in another form 
(i.e., (\ref{eq-r})), which 
is actually equivalent to (\ref{eq-r-ano}).

For any graph $G$ and any $E'\subseteq E(G)$, let 
$\Gamma(E')$ be the set of those mappings 
$g:E'\rightarrow V(G)$ 
such that for each $e\in E'$, $g(e)\in \{u_e,v_e\}$,
where $u_e$ and $v_e$ are the two ends of $e$. 
Observe that  
\begin{equation}\relabel{eq-remove-one}
\sum_{g\in \Gamma(E')}
\prod\limits_{v\in V(G)} d(v)^{-|g^{-1}(v)|}
=\prod_{e\in E'}(d(u_{e})^{-1}+d(v_{e})^{-1}).
\end{equation}

Thus (\ref{eq-r-ano}) and (\ref{eq-r0-ano})
can be replaced by the following expressions:
\begin{equation}\relabel{eq-r}
t(L(S_r(G)))=
\sum_{E'\subseteq E(G)}t(G[E'])
r^{|E'|-|V(G)|+1}
\sum_{g\in \Gamma(E(G)-E')}
\prod_{v\in V(G)} d(v)^{d(v)-2-|g^{-1}(v)|}
\end{equation}
and 
\begin{equation}\relabel{eq-r0}
t(L(G))=\sum_{T\in \sett(G)}
\sum_{g\in \Gamma(E(G)-E(T))} 
\prod_{v\in V(G)} d(v)^{d(v)-2-|g^{-1}(v)|}.
\end{equation}

\resection{Preliminary Results}

In this section, we shall establish some results which will be 
used in the next section to prove Theorem~\ref{gen-Sr} for the case $r=0$.

For any connected graph $H$ and any forest $F$ of $H$,
let $\setst_H(F)$ be the set of 
those spanning trees of $H$ containing all edges of $F$,
and $\setsf_H(F)$ be the set of those spanning forests 
of $H$ containing all edges of $F$.

In this section, we always assume that 
$G$ is a connected graph with a clique $V_0$ such that 
%containing a clique $V_0$ of order $k$ such that 
$F=G-E(G[V_0])$ is a forest and 
every vertex of $V_0$ is incident with at most one edge of $F$,
as shown in Figure~\ref{f8}.

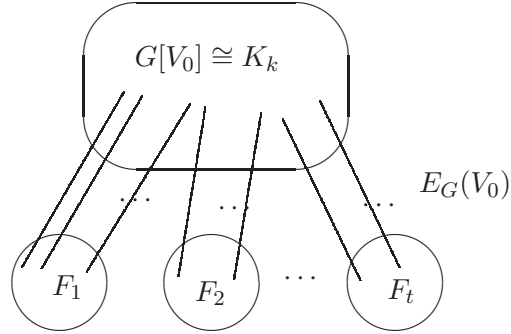
\begin{figure}[htbp]
 \centering
%%%%%%
%\input f8.pic
%%%%%%
%%%
%TeXCAD (http://texcad.sf.net/) Picture. File: [f8.pic]. Options on following lines.
%\grade{\on}
%\emlines{\off}
%\epic{\off}
%\beziermacro{\on}
%\reduce{\on}
%\snapping{\off}
%\pvinsert{% Your \input, \def, etc. here}
%\quality{8.000}
%\graddiff{0.005}
%\snapasp{1}
%\zoom{4.0000}
\unitlength .8mm % = 2.276pt
\linethickness{0.4pt}
\ifx\plotpoint\undefined\newsavebox{\plotpoint}\fi % GNUPLOT compatibility
\begin{picture}(74,59.25)(0,0)
\put(10.75,12.75){\circle{15}}
\put(36,12.75){\circle{15}}
\put(66.5,12.75){\circle{15}}
\put(48,12.25){$\cdots$}
\put(20.5,25.25){$\cdots$}
\put(61,24.5){$\cdots$}
\put(37,24){$\cdots$}
\put(23.25,48.75){$G[V_0]\cong K_k$}
\put(70.5,27.75){$E_G(V_0)$}
\put(9.5,10.75){$F_1$}
\put(33.25,9.75){$F_2$}
\put(65.25,10){$F_t$}
\put(36.75,45.375){\oval(44,27.75)[]}
%\emline(4.5,15.25)(21.5,44.5)
\multiput(4.5,15.25)(.0420792079,.0724009901){404}{\line(0,1){.0724009901}}
%\end
%\emline(7.75,15)(24.5,43.75)
\multiput(7.75,15)(.0420854271,.0722361809){398}{\line(0,1){.0722361809}}
%\end
%\emline(15.25,14.5)(32.5,42.5)
\multiput(15.25,14.5)(.0420731707,.0682926829){410}{\line(0,1){.0682926829}}
%\end
%\emline(30.5,13.75)(35,42)
\multiput(30.5,13.75)(.042056075,.264018692){107}{\line(0,1){.264018692}}
%\end
%\emline(39.75,13.25)(44.25,41)
\multiput(39.75,13.25)(.042056075,.259345794){107}{\line(0,1){.259345794}}
%\end
%\emline(60.75,13)(47.75,40)
\multiput(60.75,13)(-.0420711974,.0873786408){309}{\line(0,1){.0873786408}}
%\end
%\emline(54,43)(67.25,15.25)
\multiput(54,43)(.0420634921,-.0880952381){315}{\line(0,-1){.0880952381}}
%\end
\end{picture}
%%%
  \caption{ $V_0$ is a clique of $G$ such that 
and $G-E(G[V_0])$ is a forest}
\relabel{f8}
\end{figure}

Let $k=|V_0|$, $d=|E_G(V_0)|$, $t=c(G-V_0)$
and $F_1, F_2, \cdots, F_t$ be components of $G-V_0$.
Observe that $k\ge d\ge t$,  
as $|E_G(v, V-V_0)|\le 1$ holds for each $v\in V_0$ 
and $|E_G(V_0, V(F_i))|\ge 1$ holds for each $F_i$.

The main purpose in this section is 
to show that if $k>d$, then 
the set $\setst_G(F)$ can be equally partitioned 
into 
$\prod_{1\le j\le t} |E_G(V_0,V(F_j))|$ subsets, 
each of which has its size $k^{k-2+t-d}$.

In the following, we divide this section into two parts.

\subsection{A preliminary result on trees}

In this subsection, we shall establish some results on 
trees which are needed for the next subsection
and following sections. 

Let $T$ be any tree and $V_0$ be any proper subset of $V(T)$.
Observe that identifying all vertices in $V_0$ 
changes $T$ to a connected graph
which is a tree if and only if $|E_T(V_0)|=c(T-V_0)$.
%Now let $t=c(t-V_0)$ and $S$ be a subset of $E_T(V_0)$.
So the following observation is obvious. 

\begin{lem}\relabel{keep-t-edges}
%For any subset $S$ of $E_T(V_0)$,
Let $t=c(t-V_0)$ and $S$ be any proper subset of $E_T(V_0)$.
Then the two statements below are equivalent: 
\begin{enumerate}
\item $|S\cap E_T(V_0, V(F_i))|=1$ holds for all 
components $F_1, F_2,\cdots, F_t$ of $T-V_0$;
\item the graph obtained from $T$ by 
removing all edges in the set $E_T(V_0)-S$
and identifying all vertices of $V_0$ is a tree. 
\end{enumerate}
\end{lem}

With $T, V_0$ given above together with 
a special vertex $v\in V_0$ such that $N(v)\subseteq V_0$, 
a subset $S$ of $E_T(V_0)$
with the properties in Lemma~\ref{keep-t-edges}
will be determined by a procedure below
(i.e., Algorithm A). 
As $S$ is uniquely determined by $T, V_0$ and $v$,
 we can denote it by  $\Phi(T,V_0,v)$.
 Thus $|\Phi(T,V_0,v)|=t=c(T-V_0)$.
 
Roughly, if $t=1$, the only edge of $\Phi(T,V_0,v)$ 
will be  selected from $E_T(V_0)$ 
according to the condition that 
it has one end in the same component of $T[V_0]$ as $v$;
if 
%As $T$ is a tree, such an edge is unique. If 
$t\ge 2$, 
the $t$ edges of $\Phi(T,V_0,v)$ will be determined by the 
$t-1$ paths $P_2, P_3, \cdots,P_t$ in $T$,
where $P_j$ is the shortest path  
connecting vertices of $F_1$ and vertices of $F_j$ for $j=2,3,\cdots,t$ and 
 $F_1, F_2,\cdots, F_t$ are the components of 
$T-V_0$.

Assume  that in Algorithm A, 
$E(T)=\{e_i: i\in I\}$ for some finite 
$I$ of positive integers.

\noindent {\bf Algorithm A} with input $(T,V_0,v)$:

\begin{enumerate}
\item[Step A1.] Let $t=c(T-V_0)$. 

\item[Step A2.] If $t=1$, let $\Phi=\{e_j\}$,
where $e_j$ is the unique edge 
in the set $E_T(V_0)$
which has one end in the component of $T[V_0]$
containing $v$. Go to Step A5.

\item[Step A3.] (Now we have $t\ge 2$.)

\begin{enumerate}
\item[A3-1.]
 The components of $T-V_0$ are labeled as 
$F_1,F_2,\cdots,F_t$ such that 
\begin{equation}
\min\{s: e_s\in E_T(V_0,F_i)\}
<\min\{s': e_{s'}\in E_T(V_0,F_{i+1})\}
\end{equation}
for all $i=1,2,\cdots,t-1$.
(In other words, these components are sorted by the 
minimum edge labels. 
For example, for the tree $T$ in Figure~\ref{f11}(a),
the four components $F_1,F_2, F_3,F_4$ of $T-V_0$ 
are labeled according to this rule. )

\item[A3-2.]
For $j=2,3,\cdots,t$, %as $T$ is a tree, 
determine the unique path $P_j$ in $T$ 
which is the shortest one among all those paths in $T$ 
connecting vertices of $F_1$ to vertices of $F_j$.
\end{enumerate}

\item[Step A4.] Let 
$\Phi=
(E(P_2)\cap E_T(V_0,V(F_1)))
\cup \bigcup_{j=2}^t (E(P_j)\cap E_T(V_0,V(F_j))).
$

\item[Step A5.] Output $\Phi$.
\end{enumerate}

{\bf Remarks}: 
\begin{enumerate}
\item  Vertex $v$ is needed only for the case that $t=1$;
%\item  The vertex $v$ in $V_0$ is needed in Step A2 to determine the unique edge in $E_T(V_0,V(F_1))$ which has one end in the component of $T[V_0]$containing $v$;
\item If $t=1$, the only edge of $\Phi$ is uniquely determined as $T$ is a tree and $T-V_0$ is connected;
\item As $T$ is a tree and $F_1$ and $F_j$ are connected, 
$P_j$ is actually the only path of $T$ 
with its ends in $F_1$ and $F_j$ respectively 
and every internal vertex of $P_j$ does not belong to 
$V(F_1)\cup V(F_j)$. Thus for $P_j$ chosen in Step A3, 
$$
|E(P_j)\cap E_G(V_0,V(F_1))|
=|E(P_j)\cap E_G(V_0,V(F_j))|=1,
$$
implying that by Step A4, 
$|\Phi\cap  E_G(V_0,V(F_j))|=1$ for all $j=1,2,\cdots,t$.
\end{enumerate}

\begin{figure}[htbp]
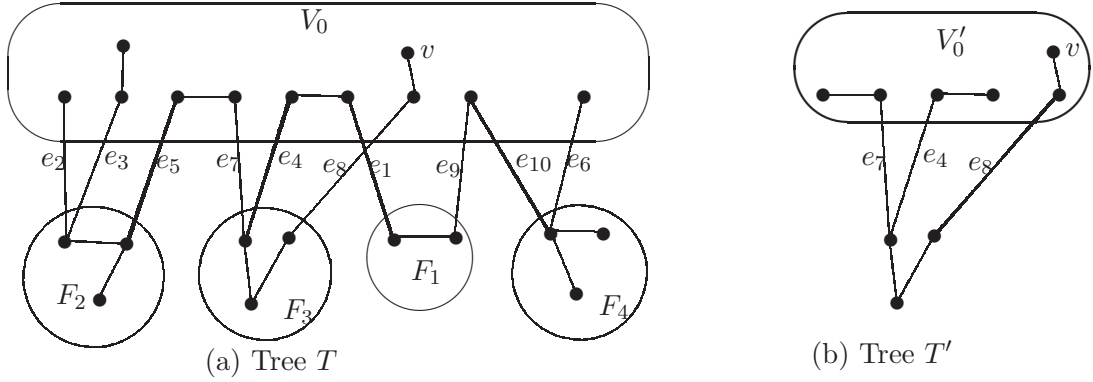

  \centering
 %\input  f11.pic

%%
%TeXCAD (http://texcad.sf.net/) Picture. File: [f11.pic]. Options on following lines.
%\grade{\on}
%\emlines{\off}
%\epic{\off}
%\beziermacro{\on}
%\reduce{\on}
%\snapping{\off}
%\pvinsert{% Your \input, \def, etc. here}
%\quality{8.000}
%\graddiff{0.005}
%\snapasp{1}
%\zoom{4.0000}
\unitlength .8mm % = 2.276pt
\linethickness{0.4pt}
\ifx\plotpoint\undefined\newsavebox{\plotpoint}\fi % GNUPLOT compatibility

%%

%%%%%%%%end of input
  \caption{$\Phi(T,V_0,v)=\{e_1,e_4,e_5,e_{10}\}$
and $\Phi(T',V'_0,v)=\{e_8\}$ }
\relabel{f11}
\end{figure}

For example, for the tree $T$ with $V_0$ and $v$ shown 
in Figure~\ref{f11} (a),
$T-V_0$ has four components $F_1,F_2,F_3,F_4$
labeled according to the minimum edge labels,
and running Algorithm A with input $(T,V_0,v)$ gives 
$\Phi(T,V_0,v)=\{e_1,e_4,e_5,e_{10}\}$, as 
the three paths $P_2,P_3$ and $P_4$ obtained by the 
algorithm 
have properties that 
$\{e_1,e_5\}\subseteq E(P_2)$,
$\{e_1,e_4\}\subseteq E(P_3)$ and
$\{e_9,e_{10}\}\subseteq E(P_4)$.
For the tree $T'$ in  Figure~\ref{f11} (b), 
$T'-V'_0$ has one component only and 
$\Phi(T',V'_0,v)=\{e_8\}$.
Note that vertex $v$ is used  for 
finding $\Phi(T',V'_0,v)$ but not for 
finding $\Phi(T,V_0,v)$.

%Note that vertex $v$ is needed only when $t=1$.

Our second purpose in this subsection is to 
show that if $|E_T(V_0,V(F_j))|>1$ for some component $F_j$ 
of $T-V_0$,
we can find another tree $T'$ with $V(T')=V(T)$ 
and $T'-E(T'[V_0])=T-E(T[V_0])$ 
such that $\Phi(T',V_0,v)$ and $\Phi(T,V_0,v)$
are different only at choosing the edge 
joining a vertex of $V_0$ to a vertex in $F_j$.

For two distinct edges $e,e'$ of $E_T(V_0)$
incident with $u$ and $u'$ respectively, 
where $u,u'\in V_0$, 
%with $u$ (resp. $u'$) as its end of $e$ (resp. of $e'$)in $V_0$,
let $T(e\leftrightarrow  e')$ be the graph, 
as shown in Figure~\ref{f12},
obtained from $T$ by 
changing every edge $(u,w)$ of $T[V_0]$, where $w\ne u'$, 
to $(u',w)$
and every edge $(u',w')$ of $T[V_0]$, where $w'\ne u$, 
to $(u,w')$.

Roughly, $T(e\leftrightarrow e')$ is actually obtained from $T$
by exchanging $(N_T(u)\cap V_0)-\{u'\}$ with 
$(N_T(u')\cap V_0)-\{u\}$.
Note that $u$ and $u'$ are adjacent in $T$ 
if and only if they are adjacent in $T(e\leftrightarrow e')$.

\begin{figure}[htbp]
  \centering
%\input f12.pic
%%
%TeXCAD (http://texcad.sf.net/) Picture. File: [f12.pic]. Options on following lines.
%\grade{\on}
%\emlines{\off}
%\epic{\off}
%\beziermacro{\on}
%\reduce{\on}
%\snapping{\off}
%\pvinsert{% Your \input, \def, etc. here}
%\quality{8.000}
%\graddiff{0.005}
%\snapasp{1}
%\zoom{4.0000}
\unitlength .8mm % = 2.276pt
\linethickness{0.4pt}
\ifx\plotpoint\undefined\newsavebox{\plotpoint}\fi % GNUPLOT compatibility
\begin{picture}(86.75,50.25)(0,0)
\put(3.75,12.75){\circle*{2.5}}
\put(63,13){\circle*{2.5}}
\put(23.25,13.25){\circle*{2.5}}
\put(82.5,13.5){\circle*{2.5}}
\put(4,33.5){\circle*{2.5}}
\put(63.25,33.75){\circle*{2.5}}
\put(23.5,34){\circle*{2.5}}
\put(82.75,34.25){\circle*{2.5}}
%\emline(3.75,33.75)(12.75,49.25)
\multiput(3.75,33.75)(.042056075,.072429907){214}{\line(0,1){.072429907}}
%\end
%\emline(3.75,34)(12,38.75)
\multiput(3.75,34)(.07300885,.042035398){113}{\line(1,0){.07300885}}
%\end
%\emline(23.25,34.5)(14,32.5)
\multiput(23.25,34.5)(-.19270833,-.04166667){48}{\line(-1,0){.19270833}}
%\end
%\emline(23.25,34.25)(14.5,23.5)
\multiput(23.25,34.25)(-.042067308,-.051682692){208}{\line(0,-1){.051682692}}
%\end
\put(12.5,48.75){\circle*{2.5}}
\put(71.75,49){\circle*{2.5}}
\put(11.75,38.5){\circle*{2.55}}
\put(71,38.75){\circle*{2.55}}
\put(14.5,32.5){\circle*{2.55}}
\put(73.75,32.75){\circle*{2.55}}
\put(14.5,23.5){\circle*{2.693}}
\put(73.75,23.75){\circle*{2.693}}
\put(.25,16.75){$e$}
\put(59.5,17){$e$}
\put(25.75,16.25){$e'$}
\put(85,16.5){$e'$}
\put(-2.75,33.75){$u$}
\put(56.5,34){$u$}
\put(27.5,34){$u'$}
\put(86.75,34.25){$u'$}
%\emline(72,49.25)(82,34.5)
\multiput(72,49.25)(.042016807,-.06197479){238}{\line(0,-1){.06197479}}
%\end
%\emline(82,34.5)(71.25,39)
\multiput(82,34.5)(-.10046729,.042056075){107}{\line(-1,0){.10046729}}
%\end
%\emline(73.5,33)(62.75,34)
\multiput(73.5,33)(-.4479167,.0416667){24}{\line(-1,0){.4479167}}
%\end
%\emline(63.25,33.75)(73.75,23.75)
\multiput(63.25,33.75)(.044117647,-.042016807){238}{\line(1,0){.044117647}}
%\end
\put(6,1){Part of $T$}
\put(62.75,3.25){Part of $T(e\leftrightarrow e')$}
\put(13.75,44.5){\makebox(0,0)[cc]{$\vdots$}}
\put(72.25,44.5){\makebox(0,0)[cc]{$\vdots$}}
\put(14.5,29.25){\makebox(0,0)[cc]{$\vdots$}}
\put(73,29.25){\makebox(0,0)[cc]{$\vdots$}}
\put(3.5,12.75){\line(0,1){20.75}}
\put(23.25,13.25){\line(0,1){21}}
\put(63,13.25){\line(0,1){21}}
\put(82.25,13.75){\line(0,1){20.75}}
\end{picture}

  \caption{$T$ and $T(e\leftrightarrow e')$}
\relabel{f12}
\end{figure}
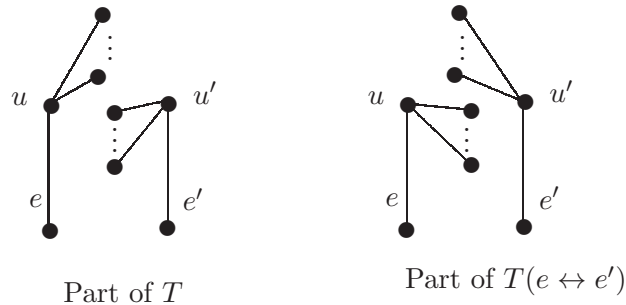

Let $T'$ denote $T(e\leftrightarrow e')$ in the remainder of 
this subsection.
There is a bijection $\tau: E(T)\rightarrow E(T')$
 defined below: 
$\tau(e)=e'$, $\tau(e')=e$,
$\tau((u,w))=(u',w)$ whenever 
$(u,w)\in E(T)$ for $w\in V_0-\{u'\}$,
$\tau((u',w'))=(u,w')$ whenever $(u',w')\in E(T)$ 
for $w'\in V_0-\{u\}$,
and $\tau(e'')=e''$ for all other edges $e''$ in $T$.

Note that $T'$ may be not a tree,
although $T'-V_0$ and $T-V_0$ are the same graph
and $F_1,F_2,\cdots,F_t$ are 
also the components of $T'-V_0$.
But $T'$ is indeed a tree if both $e$ and $e'$ 
have ends in the same component of $T-V_0$.

\begin{lem}\relabel{comp-twotrees}
Let $e,e'$ be distinct edges of $E_T(V_0,V(F_i))$ for some 
$i$ with $1\le i\le t$.
%let $T'$ denote $T(e\leftrightarrow e')$.
\begin{enumerate}
\item Then $T'$ is a tree;
\item If $e\in \Phi(T,V_0,v)$ and 
either $t\ge 2$ or $N_T(v)\subseteq V_0$,
then 
$\Phi(T',V_0,v)=(\Phi(T,V_0,v)-\{e\})\cup \{e'\}$.
\end{enumerate}
\end{lem}

\proof Note that  for any edge $e''\in E(T-V_0)$,
$T/e''$ is also a tree, 
$T'$ is a tree if and $T'/e''$ is a tree, and 
more importantly, 
$\Phi(T,V_0,v)=\Phi(T/e'',V_0,v)$.
Thus it suffices to prove this lemma only for 
the case that $|V(F_i)|=1$ for all $i=1,2,\cdots,t$. 

%Let $T'$ denote $T(e\leftrightarrow e')$.

(i) It can be proved easily by induction on the number of edges 
in $T$.

(ii) Assume that $t=1$. 
Then $N_T(v)\subseteq V_0$ and so $v$ is not any end of $e$.
As $e\in  \Phi(T,V_0,v)$,
$\Phi(T,V_0,v)=\{e\}$. By Algorithm A, 
$e$ has one end (i.e., $u$) in the component of $T[V_0]$ 
containing $v$
(i.e., 
the subgraph $T[V_0]$ has a path $P$ 
connecting $v$ to $u$).
By the definition of $T'$ (i.e., $T(e\leftrightarrow e')$),
$P$ is now changed to a path $P'$ in 
$T'[V_0]$ by the mapping $\tau$ 
connecting $v$ to one end of $e'$ (i.e., $u'$).
Thus %$e'$ has one end in the component of $T'[V_0]$ containing $v$.
 $\Phi(T',V_0,v)=\{e'\}$ by Algorithm A.
The result holds for this case.

Now assume that $t\ge 2$.
For $j=2,3,\cdots,t$, 
let $P_j$ be the only path in $T$ 
with its two ends in $F_1$ and $F_j$ respectively 
and every interval vertex of $P_j$ 
does not below to $V(F_1)\cup V(F_j)$.

With the bijection $\tau: E(T)\rightarrow E(T')$
defined above, for $j=2,3,\cdots,t$,
$\tau(E(P_j))$ is a subset of $E(T')$ 
and forms a path in $T'$, denoted by $P'_j$.
Note that the two ends of $P'_j$ are in 
$F_1$ and $F_j$ respectively 
and every interval vertex of $P'_j$ 
does not below to $V(F_1)\cup V(F_j)$.
Also observe that for $j=2,3,\cdots,t$,
if $i\in \{1,j\}$, then
$$
E(P'_j)\cap E_{T'}(V_0,V(F_i))=\{e'\},
$$
and if $s\in \{1,j\}-\{i\}$, then 
$$
E(P'_j)\cap E_{T'}(V_0,V(F_s))=
E(P_j)\cap E_{T}(V_0,V(F_s)).
$$
Hence (ii) holds. 
\proofend

\subsection{Partitions of $\setst_G(F)$}

Recall that $G$ is a connected graph 
with a clique $V_0$ of order $k$ such that $F=G-E(G[V_0])$ is a forest and every vertex of $V_0$ is incident with at most one edge of $F$
(i.e., $d_F(v)\le 1$ for each $v\in V_0$),
as shown in Figure~\ref{f8}. 
In this subsection, 
our main purpose is to  partition 
$\setst_G(F)$ equally into 
$\prod_{j=1}^t |E_G(V_0,V(F_j))|$ subsets,
where $F_1,F_2,\cdots,F_t$ are the components of $G-V_0$.

We start with the following beautiful formula for
the number of spanning trees of 
a complete graph $K_k$ of order $k$ containing a given 
spanning forest.
This result was originally 
due to Lov\'{a}sz  (Problem 4 in page 29 of \cite{lov}).
%This result can also be found in (\cite{Mey2011}).

\begin{theo}[\cite{lov}]\relabel{forest-tree}
For any spanning forest $F$ of $K_k$, 
if $c$ is the number of components of $F$ 
and $k_1,k_2,\cdots,k_c$ are the orders of its components, then 
$$
|\setst_{K_k}(F)|=k^{c-2}\prod_{i=1}^c k_i.
$$
\end{theo}

This result naturally generalizes the well-known formula that 
$t(K_k)=k^{k-2}$ for any $k\ge 1$, 
which was first obtained by Cayley~\cite{Aig}.
%Recall that in this section, $G$ is a graph containing a clique $V_0$ of order $k$ such that $F=G-E(G[V_0])$ is a forest and  every vertex of $V_0$ is incident with at most one edge of $F$,as shown in Figure~\ref{f8}.
Now we apply this result %is applied in this subsection %Theorem~\ref{forest-tree} 
to establish some results on
the set $\setst_G(F)$
and finally partition 
$\setst_G(F)$ equally into 
$\prod_{j=1}^t |E_G(V_0,V(F_j))|$ subsets.

Recall that $d=|E_G(V_0)|$ and $k\ge d\ge t$.

\begin{pro}\relabel{forest-stree} 
With $G, F$ and $V_0$ defined above,  we have 
$$
|\setst_G(F)|=k^{k-2+t-d}
\prod_{j=1}^t |E_G(V_0,V(F_j))|.
$$
%where the product runs over all $t$ components $F_j$ of $G-V_0$.
\end{pro}

\proof 
We only need to consider the case that 
$E_G(V_0,V(F_j))\ne \emptyset$ for every component $F_j$ 
of $G-V_0$;
otherwise, the result is trivial 
as $|\setst_G(F)|=0$ when $G$ is disconnected.

Observe that for any edge $e$ of 
$E(G-V_0)$, 
we have $|\setst_{G/e}(F/e)|=|\setst_G(F)|$.
Thus we may assume that every component of $G-V_0$ 
is a single vertex, implying that 
$G-V_0$ is the empty graph of $t$ vertices,
namely $x_1,x_2,\cdots,x_t$.
So $E(F)=E_G(V_0)$.

For each $j=1,2,\cdots,t$, let 
$V_j=\{x\in V_0: x \mbox{ is incident with }x_j\}$
and $e_j$ be any edge joining $x_j$ to some vertex in $V_j$.
Let $G_0=G[V_0]$. Note that 
$F/\{e_1,e_2,\cdots,e_t\}$ 
can be considered as a spanning forest of $G_0$ 
and 
$$
\setst_G(F)=\setst_{G_0}(F/\{e_1,e_2,\cdots,e_t\}).
$$
As $G_0$ is a complete graph of order $k$, 
by Theorem~\ref{forest-tree}, 
$$
|\setst_{G_0}(F/\{e_1,e_2,\cdots,e_t\})|
=k^{c-2}\prod_{j=1}^c |V'_j|,
$$
where $c$ is the number of components of $F/\{e_1,e_2,\cdots,e_t\}$
and $V'_1, V'_2, \cdots, V'_c$ are vertex sets of components
of  $F/\{e_1,e_2,\cdots,e_t\}$.
Note that 
%There are exactly $k-d$ vertices in the following set 
$$
|V_0-\bigcup_{j=1}^t V_j|
=|V_0|-\sum_{k=1}^t |V_j|
=k-|E_G(V_0)|
=k-d,
$$
implying that  $c=k-d+t$
and the sizes of $V_1', V_2', \cdots, V_c'$ are 
equal to 
$$
|V_1|,\cdots, |V_t|, \underbrace{1,\cdots,1}_{k-d}.
$$
As $|V_j|=|E_G(V_0,\{x_j\})|=|E_G(V_0,V(F_j))|$,  
the result follows from Theorem~\ref{forest-tree}. 
\proofend

Now assume that %Assume that $G$ is connected,
$v$ is a vertex of $V_0$
with $N_G(v)\subseteq V_0$,
i.e., $d_F(v)=0$.
Note that this condition %$N_G(v)\subseteq V_0$
is only needed for the case that $G-V_0$ is connected. 
Under this condition, it is obvious that $k>d$.

Recall that for any tree $T$ of $\setst_G(F)$, 
$\Phi(T,V_0,v)$ is a subset of $E_G(V_0)$ 
with the property that
$|\Phi(T,V_0,v)\cap E_G(V_0,V(F_j))|=1$ 
for each $j=1,2,\cdots,t$.
For each subset $S$ of $E_G(V_0)$ with 
the property that 
$|S\cap E_G(V_0,V(F_j))|=1$ for each $j=1,2,\cdots,t$,
let $\setst_G(F,S,v)$ denote the set of those 
spanning trees $T\in \setst_G(F)$ with 
$\Phi(T,V_0,v)=S$. 
Thus $\setst_G(F)$ is partitioned into 
$\prod_{j=1}^t |E_G(V_0,V(F_j))|$ subsets
$\setst_G(F,S,v)$'s.
The following result shows that 
all these sets $\setst_G(F,S,v)$'s
have the same size.

The following result shows that $|\setst_G(F,S,v)|$
is independent of $S$.

\begin{pro}\relabel{no-tree-ext}
Assume that $k>d$ and $N(v)\subseteq V_0$.
For any subset $S$ of $E_G(V_0)$
with $|S\cap E_G(V_0,V(F_j))|=1$ 
for each component $F_j$ of $G-V_0$,
%For any set $S$ of $t$ edges, one from $E_G(V_0,V(F_j))$ for each component $F_j$ of $G-V_0$,
we have 
$$
|\setst_G(F,S,v)|=k^{k-2+t-d}.
$$
\end{pro}

\proof There are exactly
$\prod_{j=1}^t |E_G(V_0,V(F_j))|$ subsets $S$ of $E_G(V_0)$
with the property that 
$|S\cap E_G(V_0,V(F_j))|=1$ for each component $F_j$ of $G-V_0$.
By Proposition~\ref{forest-stree}, 
we only need to show that 
$|\setst_G(F,S,v)|=|\setst_G(F,S',v)|$ holds for any two such sets $S$ and $S'$.
Thus it suffices to show that $|\setst_G(F,S,v)|=|\setst_G(F,S',v)|$ holds 
for any two such sets $S$ and $S'$ with $|S-S'|=1$,
i.e., $S$ and $S'$ have exactly $t-1$ same edges.

Let $S$ be such a subset of $E_G(V_0)$ mentioned above. 
Assume that $e,e'$ are distinct edges in 
$E_G(V_0,V(F_j))$ for some $j$ 
with $e\in S$ and $e'\notin S$.
Let $S'=(S-\{e\})\cup \{e'\}$.
It remains to show that 
$|\setst_G(F,S,v)|=|\setst_G(F,S',v)|$.

For any $T\in \setst_G(F,S,v)$, 
let $T'$ be the tree $T(e\leftrightarrow e')$.
By Lemma~\ref{comp-twotrees}, 
we have $\Phi(T',V_0,v)=(\Phi(T,V_0,v)-\{e\})\cup \{e'\}$,
implying that $T'\in \setst_G(F,S',v)$.

Let $\phi$ be the mapping from $\setst_G(F,S,v)$ to  $\setst_G(F,S',v)$ defined by $\phi(T)=T(e\leftrightarrow e')$.
It is obvious that $\phi$ is an onto mapping,
and $\phi': T' \rightarrow T'(e'\leftrightarrow e)$ 
is also an onto mapping 
from $\setst_G(F,S',v)$ to  $\setst_G(F,S,v)$.

Thus $|\setst_G(F,S,v)|=|\setst_G(F,S',v)|$ 
and  the result follows.
\proofend

We end this section with an application of 
Proposition~\ref{no-tree-ext} to deduce another result.

Let $G'$ be the graph obtained from $G$ by contracting 
all edges in $G[V_0]$.
Then $V_0$ becomes a vertex in $G'$,
denoted by $v_0$.  
Thus $V(G')=(V(G)-V_0)\cup \{v_0\}$,
and $E(G')$ and $E(G)-E(G[V_0])$ are the same 
although for each edge $e\in E_G(V_0)$,
its end in $V_0$ is changed to $v_0$ 
when $e$ becomes an edge in $G'$.
An example for $G$ and $G'$ 
is shown in Figure~\ref{f14} (a) and (b).

Let $T'$ be any spanning tree of $G'$ 
with $E(G'-v_0)\subseteq E(T')$, i.e.,  
$T'\in \setst_{G'}(F')$ for $F'=G'-v_0$. 
Thus $|E_{T'}(v_0)|=t$, i.e., $E_{T'}(v_0)$
has exactly $t$ edges, 
corresponding to $t$ edges in $G$, 
one from $E_G(V_0,V(F_j))$ for each component $F_j$ of $G-V_0$.
An example for $T'$ is shown in Figure~\ref{f14} (d).

Let $D$ be any subset of  $E_G(V_0)-E_{T'}(v_0)$
and 
let $\setst_G(V_0,T',D,v)$ be 
the set of those spanning trees $T$ of $G$ 
such that 
(i) $T-V_0$ and $T'-v_0$ are the same graph;
(ii) $E_T(V_0)$ is the disjoint union of  $D$ and $E_{T'}(v_0)$
and 
(iii) $\Phi(T,V_0,v)=E_{T'}(v_0)$.
Thus $\setst_G(V_0,T',D,v)\subseteq \setst_G(F)$ if and only if $D=E_G(V_0)-E_{T'}(v_0)$.
For example, the tree $T$ in Figure~\ref{f14} (c)
belongs to $\setst_G(V_0,T',D,v)$ with $D=\{e_1,e_5\}$,
but $T\notin \setst_G(F)$, 
as $E(F)\not\subseteq E(T)$.

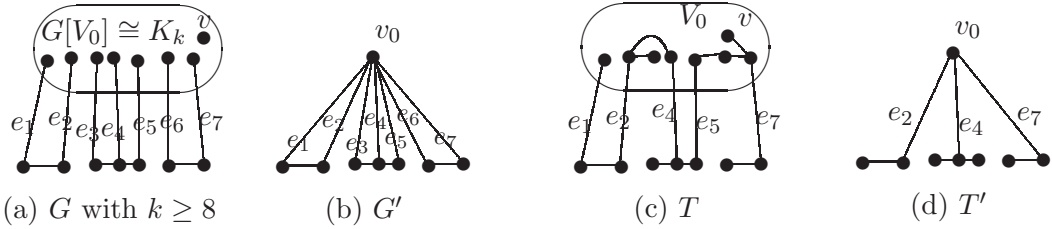
\begin{figure}[htbp]
  \centering
 %\input  f14.pic

%%
%TeXCAD (http://texcad.sf.net/) Picture. File: [f14.pic]. Options on following lines.
%\grade{\on}
%\emlines{\off}
%\epic{\off}
%\beziermacro{\on}
%\reduce{\on}
%\snapping{\off}
%\pvinsert{% Your \input, \def, etc. here}
%\quality{8.000}
%\graddiff{0.005}
%\snapasp{1}
%\zoom{6.7272}
\unitlength .8mm % = 2.276pt
\linethickness{0.4pt}
\ifx\plotpoint\undefined\newsavebox{\plotpoint}\fi % GNUPLOT compatibility
\begin{picture}(172.531,35.548)(0,0)
\put(19.375,28){\oval(31.25,14.5)[]}
\put(110.349,28.298){\oval(31.25,14.5)[]}
%\emline(5.75,26.25)(2,8.5)
\multiput(5.75,26.25)(-.04213483,-.1994382){89}{\line(0,-1){.1994382}}
%\end
%\emline(98.365,26.399)(94.615,8.649)
\multiput(98.365,26.399)(-.04213483,-.1994382){89}{\line(0,-1){.1994382}}
%\end
\put(2,8.5){\line(1,0){6.5}}
\put(94.615,8.649){\line(1,0){6.5}}
\put(45.108,8.554){\line(1,0){6.5}}
\put(141.5,9.25){\line(1,0){6.5}}
%\emline(8.5,8.5)(10,25.75)
\multiput(8.5,8.5)(.04166667,.47916667){36}{\line(0,1){.47916667}}
%\end
%\emline(101.115,8.649)(102.615,25.899)
\multiput(101.115,8.649)(.04166667,.47916667){36}{\line(0,1){.47916667}}
%\end
\put(13.75,8.75){\line(1,0){7.5}}
\put(106.365,8.899){\line(1,0){7.5}}
\put(56.858,8.804){\line(1,0){7.5}}
\put(153.25,9.5){\line(1,0){7.5}}
\put(26.25,8.75){\line(1,0){5}}
\put(118.865,8.899){\line(1,0){5}}
\put(69.358,8.804){\line(1,0){5}}
\put(165.75,9.5){\line(1,0){5}}
\put(6,26){\circle*{2.062}}
\put(98.615,26.149){\circle*{2.062}}
\put(2,8.25){\circle*{2.062}}
\put(94.615,8.399){\circle*{2.062}}
\put(45.108,8.304){\circle*{2.062}}
\put(141.5,9){\circle*{2.062}}
\put(8.75,8.25){\circle*{2.062}}
\put(101.365,8.399){\circle*{2.062}}
\put(51.858,8.304){\circle*{2.062}}
\put(148.25,9){\circle*{2.062}}
\put(14,8.75){\circle*{2.062}}
\put(106.615,8.899){\circle*{2.062}}
\put(57.108,8.804){\circle*{2.062}}
\put(153.5,9.5){\circle*{2.062}}
\put(18,8.75){\circle*{2.062}}
\put(110.615,8.899){\circle*{2.062}}
\put(61.108,8.804){\circle*{2.062}}
\put(157.5,9.5){\circle*{2.062}}
\put(21.25,8.75){\circle*{2.062}}
\put(113.865,8.899){\circle*{2.062}}
\put(64.358,8.804){\circle*{2.062}}
\put(160.75,9.5){\circle*{2.062}}
\put(26.25,8.5){\circle*{2.062}}
\put(118.865,8.649){\circle*{2.062}}
\put(69.358,8.554){\circle*{2.062}}
\put(165.75,9.25){\circle*{2.062}}
\put(32,8.75){\circle*{2.062}}
\put(124.615,8.899){\circle*{2.062}}
\put(75.108,8.804){\circle*{2.062}}
\put(171.5,9.5){\circle*{2.062}}
\put(30.25,26.25){\circle*{2.062}}
\put(122.865,26.399){\circle*{2.062}}
\put(26,26.5){\circle*{2.062}}
\put(118.615,26.649){\circle*{2.062}}
\put(21,26){\circle*{2.062}}
\put(113.615,26.149){\circle*{2.062}}
\put(17,26.5){\circle*{2.062}}
\put(109.615,26.649){\circle*{2.062}}
\put(60.108,26.554){\circle*{2.062}}
\put(156.5,27.25){\circle*{2.062}}
\put(14.25,26.5){\circle*{2.062}}
\put(106.865,26.649){\circle*{2.062}}
\put(10,26.5){\circle*{2.062}}
\put(102.615,26.649){\circle*{2.062}}
%\emline(13.5,8.75)(14.25,27)
\multiput(13.5,8.75)(.0416667,1.0138889){18}{\line(0,1){1.0138889}}
%\end
%\emline(17.25,27)(18,8.75)
\multiput(17.25,27)(.0416667,-1.0138889){18}{\line(0,-1){1.0138889}}
%\end
%\emline(109.865,27.149)(110.615,8.899)
\multiput(109.865,27.149)(.0416667,-1.0138889){18}{\line(0,-1){1.0138889}}
%\end
%\emline(60.358,27.054)(61.108,8.804)
\multiput(60.358,27.054)(.0416667,-1.0138889){18}{\line(0,-1){1.0138889}}
%\end
%\emline(156.75,27.75)(157.5,9.5)
\multiput(156.75,27.75)(.0416667,-1.0138889){18}{\line(0,-1){1.0138889}}
%\end
%\emline(21.25,8.75)(21,26)
\multiput(21.25,8.75)(-.041667,2.875){6}{\line(0,1){2.875}}
%\end
%\emline(113.865,8.899)(113.615,26.149)
\multiput(113.865,8.899)(-.041667,2.875){6}{\line(0,1){2.875}}
%\end
\put(26,8.75){\line(0,1){18.25}}
%\emline(30.25,26.75)(32,9)
\multiput(30.25,26.75)(.04166667,-.42261905){42}{\line(0,-1){.42261905}}
%\end
%\emline(122.865,26.899)(124.615,9.149)
\multiput(122.865,26.899)(.04166667,-.42261905){42}{\line(0,-1){.42261905}}
%\end
\put(17,30.75){\makebox(0,0)[cc]{$G[V_0]\cong K_k$}}
\put(16.25,1){\makebox(0,0)[cc]{(a) $G$ with $k\ge 8$}}
\put(108.865,1.149){\makebox(0,0)[cc]{(c) $T$}}
\put(59.358,1.054){\makebox(0,0)[cc]{(b) $G'$ }}
\put(155.75,1.75){\makebox(0,0)[cc]{(d) $T'$}}
\put(2,15.25){\makebox(0,0)[cc]{$e_1$}}
\put(94.615,15.399){\makebox(0,0)[cc]{$e_1$}}
\put(47.784,12.331){\makebox(0,0)[cc]{{\small $e_1$}}}
\put(8.25,16){\makebox(0,0)[cc]{$e_2$}}
\put(100.865,16.149){\makebox(0,0)[cc]{$e_2$}}
\put(53.588,15.608){\makebox(0,0)[cc]{{\small $e_2$}}}
\put(147.75,16.75){\makebox(0,0)[cc]{$e_2$}}
\put(12.75,14){\makebox(0,0)[cc]{$e_3$}}
\put(57.493,11.824){\makebox(0,0)[cc]{{\small $e_3$}}}
\put(17,14.25){\makebox(0,0)[cc]{$e_4$}}
\put(60.554,15.939){\makebox(0,0)[cc]{{\small $e_4$}}}
\put(159.5,15){\makebox(0,0)[cc]{$e_4$}}
\put(22.25,15.25){\makebox(0,0)[cc]{$e_5$}}
\put(115.865,15.399){\makebox(0,0)[cc]{$e_5$}}
\put(64.02,12.777){\makebox(0,0)[cc]{{\small $e_5$}}}
\put(26.75,15.5){\makebox(0,0)[cc]{$e_6$}}
\put(65.993,16.595){\makebox(0,0)[cc]{{\small $e_6$}}}
\put(33.25,16){\makebox(0,0)[cc]{$e_7$}}
\put(125.865,16.149){\makebox(0,0)[cc]{$e_7$}}
\put(72.196,12.784){\makebox(0,0)[cc]{{\small $e_7$}}}
\put(169.25,16.75){\makebox(0,0)[cc]{$e_7$}}
%\emline(107.115,26.899)(102.365,26.649)
\multiput(107.115,26.899)(-.791667,-.041667){6}{\line(-1,0){.791667}}
%\end
%\emline(118.615,27.149)(113.115,26.649)
\multiput(118.615,27.149)(-.458333,-.041667){12}{\line(-1,0){.458333}}
%\end
%\emline(122.865,26.649)(118.865,26.899)
\multiput(122.865,26.649)(-.666667,.041667){6}{\line(-1,0){.666667}}
%\end
%\emline(44.858,8.804)(60.358,27.054)
\multiput(44.858,8.804)(.0421195652,.0495923913){368}{\line(0,1){.0495923913}}
%\end
%\emline(156.75,27.75)(171.5,9.5)
\multiput(156.75,27.75)(.0421428571,-.0521428571){350}{\line(0,-1){.0521428571}}
%\end
%\emline(69.108,8.554)(59.858,26.804)
\multiput(69.108,8.554)(-.042045455,.082954545){220}{\line(0,1){.082954545}}
%\end
%\emline(59.858,26.804)(51.608,8.554)
\multiput(59.858,26.804)(-.042091837,-.093112245){196}{\line(0,-1){.093112245}}
%\end
%\emline(156.25,27.5)(148,9.25)
\multiput(156.25,27.5)(-.042091837,-.093112245){196}{\line(0,-1){.093112245}}
%\end
%\emline(57.108,9.054)(60.108,27.304)
\multiput(57.108,9.054)(.04166667,.25347222){72}{\line(0,1){.25347222}}
%\end
\put(60.358,29.804){$v_0$}
\put(156.75,30.5){$v_0$}
%\emline(64.358,9.304)(60.358,26.804)
\multiput(64.358,9.304)(-.042105263,.184210526){95}{\line(0,1){.184210526}}
%\end
\put(106.365,16.649){$e_4$}
\qbezier(109.365,26.649)(106.74,33.524)(102.615,26.899)
\put(110.893,32.109){$V_0$}
\put(31.96,29.879){\circle*{1.994}}
\put(119.069,30.176){\circle*{1.994}}
\put(30.771,31.514){$v$}
\put(120.704,31.812){$v$}
%\emline(122.636,26.758)(118.92,30.622)
\multiput(122.636,26.758)(-.04175281,.04341573){89}{\line(0,1){.04341573}}
%\end
%\emline(60.055,26.757)(75.217,8.919)
\multiput(60.055,26.757)(.0421175724,-.0495500852){360}{\line(0,-1){.0495500852}}
%\end
\end{picture}
%%

%%%%%%%%end of input
  \caption{A tree $T$ in $\setst_G(V_0,T',D,v)$
with $D=\{e_1,e_5\}$}
%{$T$ (resp. $T'$) is a spanning tree of $G$(resp. $G'$) and $D=E_{T}(V_0)-E_{T'}(v_0)$}
\relabel{f14}
\end{figure}

\begin{pro}\relabel{no-trees-original}
With $T'$ and $D$ given above, 
we have 
$$
|\setst_G(V_0,T',D,v)|=k^{k-2-|D|}.
$$
\end{pro}

\proof Let $G^*$ denote the graph $G-D'$,
where $D'=E_G(V_0)-(D\cup E_{T'}(v_0))$.
Observe that 
$$
\setst_G(V_0,T',D,v)
=\setst_{G^*}(F^*,E_{T'}(v_0), v),
$$
where $F^*=G^*-E(G^*[V_0])$, i.e., $F^*=F-D'$.
Also note that $c(G^*-V_0)=c(G-V_0)=t$ and 
$$
|E_{G^*}(V_0)|=|E_{T'}(v_0)\cup D|=t+|D|.
$$ 
By Proposition~\ref{no-tree-ext}, we have 
$$
|\setst_{G^*}(F^*,E_{T'}(v_0), v)|
=k^{k-2+t-(t+|D|)}
=k^{k-2-|D|}.
$$
Thus the result holds.
\proofend

\resection{Proving Theorem~\ref{gen-Sr} for $r=0$}

In this section, we shall prove
Theorem~\ref{gen-Sr} for the case $r=0$
(i.e., the result of (\ref{eq-r0}) or 
equivalently (\ref{eq-r0-ano})) 
is a special case of another result
(i.e., Theorem~\ref{gen-result}). % than that. 

Let $u$ be any vertex in a simple graph $G$.
Assume that $E_G(u)=\{(u,u_i): 1\le i\le s\}$,
where $s=d_G(u)$.
If $G'$ is the graph obtained from $G-u$ by 
adding a complete graph $K_s$ 
with vertices $w_1,w_2,\cdots,w_s$  
and adding $s$ new 
edges $(w_i,u_i)$ for $i=1,2,\cdots,s$,
then $G'$ is said to be obtained from $G$ by a 
{\it clique-insertion at $u$}.
The clique-insertion is a graph operation 
playing an important role in the study of vertex-transitive graphs (see \cite{lov2,mad}).
The {\it clique-inserted graph} of $G$, denoted by $C(G)$,
is obtained from $G$ by operating clique-insertion at 
every vertex of $G$.
Note that the clique-inserted graph of $G$ is also called 
{\it the para-line graph} of $G$ (see \cite{shi}).
An example for $C(G)$ is shown in Figure~\ref{f13}.

Let $M$ be the set of those edges in $E(C(G))$
which are not in the inserted cliques.
So $M$ consists of all edges in $E(G)$ 
and thus can be considered as the same as $E(G)$.
Observe that $C(G)$ has the following properties:
\begin{enumerate}
%\item The vertex set of $C(G)$ is$\{v_1(e),v_2(e): e\in E(G)\}$; 
\item $M$ is a matching of $C(G)$; 

\item $L(G)$ is the graph $C(G)/M$
and thus $t(L(G))=|\setst_{C(G)}(M)|$; 

\item each component of $C(G)-M$
is a complete graph. 
\end{enumerate}

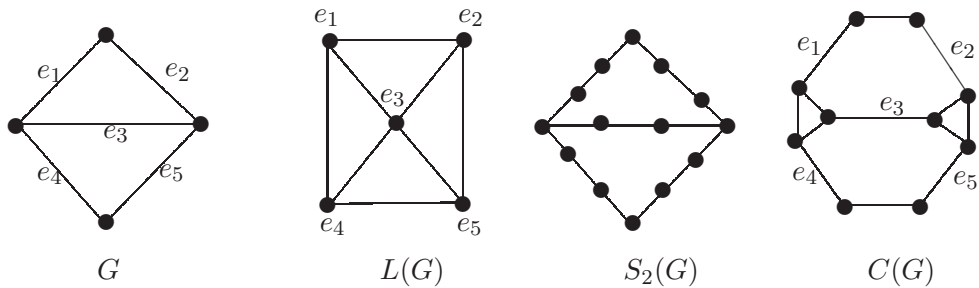
\begin{figure}[htbp]
  \centering
 \unitlength .8mm % = 2.276pt
\linethickness{0.4pt}
\ifx\plotpoint\undefined\newsavebox{\plotpoint}\fi % GNUPLOT compatibility
\begin{picture}(161.096,46.596)(0,0)
%\emline(16.25,42.5)(1.75,27.75)
\multiput(16.25,42.5)(-.0421511628,-.042877907){344}{\line(0,-1){.042877907}}
%\end
%\emline(103.75,42.25)(89.25,27.5)
\multiput(103.75,42.25)(-.0421511628,-.042877907){344}{\line(0,-1){.042877907}}
%\end
%\emline(1.75,27.75)(16.25,11.25)
\multiput(1.75,27.75)(.0421511628,-.0479651163){344}{\line(0,-1){.0479651163}}
%\end
%\emline(89.25,27.5)(103.75,11)
\multiput(89.25,27.5)(.0421511628,-.0479651163){344}{\line(0,-1){.0479651163}}
%\end
%\emline(16.25,11.25)(32.25,27.75)
\multiput(16.25,11.25)(.0421052632,.0434210526){380}{\line(0,1){.0434210526}}
%\end
%\emline(103.75,11)(119.75,27.5)
\multiput(103.75,11)(.0421052632,.0434210526){380}{\line(0,1){.0434210526}}
%\end
%\emline(32.25,27.75)(16.5,42.5)
\multiput(32.25,27.75)(-.045,.0421428571){350}{\line(-1,0){.045}}
%\end
%\emline(119.75,27.5)(104,42.25)
\multiput(119.75,27.5)(-.045,.0421428571){350}{\line(-1,0){.045}}
%\end
\put(53.25,41.25){\line(0,-1){26.75}}
%\emline(53.25,14.5)(75.75,14.75)
\multiput(53.25,14.5)(3.75,.041667){6}{\line(1,0){3.75}}
%\end
\put(75.75,14.75){\line(0,1){27}}
\put(75.75,41.75){\line(-1,0){23}}
%\emline(52.75,41.75)(75.75,14.75)
\multiput(52.75,41.75)(.0421245421,-.0494505495){546}{\line(0,-1){.0494505495}}
%\end
%\emline(75.75,41.75)(53.75,14.5)
\multiput(75.75,41.75)(-.0421455939,-.0522030651){522}{\line(0,-1){.0522030651}}
%\end
\put(141,45.5){\line(1,0){10.25}}
\put(151.25,45.5){\line(2,-3){8.5}}
\put(159.75,32.75){\line(0,-1){8.25}}
%\emline(159.75,24.5)(151.75,14.25)
\multiput(159.75,24.5)(-.042105263,-.053947368){190}{\line(0,-1){.053947368}}
%\end
\put(151.75,14.25){\line(-1,0){12.75}}
%\emline(139,14.25)(131.5,24.75)
\multiput(139,14.25)(-.042134831,.058988764){178}{\line(0,1){.058988764}}
%\end
\put(131.5,24.75){\line(0,1){8.75}}
%\emline(131.5,33.5)(140.75,45.25)
\multiput(131.5,33.5)(.042045455,.053409091){220}{\line(0,1){.053409091}}
%\end
%\emline(131.75,33.5)(136.75,28.75)
\multiput(131.75,33.5)(.044247788,-.042035398){113}{\line(1,0){.044247788}}
%\end
%\emline(136.75,28.75)(131.5,24.75)
\multiput(136.75,28.75)(-.055263158,-.042105263){95}{\line(-1,0){.055263158}}
%\end
\put(136.5,28.75){\line(1,0){17}}
%\emline(153.5,28.75)(159.25,24.75)
\multiput(153.5,28.75)(.060526316,-.042105263){95}{\line(1,0){.060526316}}
%\end
%\emline(153.5,28.5)(159.5,32.75)
\multiput(153.5,28.5)(.059405941,.042079208){101}{\line(1,0){.059405941}}
%\end
\put(2,27.75){\line(1,0){30.25}}
\put(89.5,27.5){\line(1,0){30.25}}
\put(16.5,42.5){\circle*{2.693}}
\put(104,42.25){\circle*{2.693}}
\put(99,37.5){\circle*{2.693}}
\put(95,32.75){\circle*{2.693}}
\put(98.75,28){\circle*{2.693}}
\put(108.75,27.5){\circle*{2.693}}
\put(108.75,37.5){\circle*{2.693}}
\put(115.5,31.75){\circle*{2.693}}
\put(114.5,21.75){\circle*{2.693}}
\put(109,16.75){\circle*{2.693}}
\put(98.75,16.75){\circle*{2.693}}
\put(93.25,22.75){\circle*{2.693}}
\put(32.25,27.75){\circle*{2.693}}
\put(119.75,27.5){\circle*{2.693}}
\put(16.5,11.5){\circle*{2.693}}
\put(104,11.25){\circle*{2.693}}
\put(1.5,27.5){\circle*{2.693}}
\put(89,27.25){\circle*{2.693}}
\put(53.75,41.5){\circle*{2.693}}
\put(76,41.75){\circle*{2.693}}
\put(64.75,28){\circle*{2.693}}
\put(53.25,14.25){\circle*{2.693}}
\put(75.75,14.5){\circle*{2.693}}
\put(131,25){\circle*{2.693}}
\put(131.75,33.75){\circle*{2.693}}
\put(136.5,29){\circle*{2.693}}
\put(154.25,28.5){\circle*{2.693}}
\put(159.75,32.5){\circle*{2.693}}
\put(159.75,24){\circle*{2.693}}
\put(151.75,14){\circle*{2.693}}
\put(139.25,14){\circle*{2.693}}
\put(141.25,45.25){\circle*{2.693}}
\put(151.25,45){\circle*{2.693}}
\put(15,2){$G$}
\put(102.5,2){$S_2(G)$}
\put(62,2){$L(G)$}
\put(5.00,35.5){$e_1$}
\put(26.25,35.75){$e_2$}
\put(16,25.5){$e_3$}
\put(5.00,18.5){$e_4$}
\put(25.25,19){$e_5$}
\put(51,44.75){$e_1$}
\put(75,44.75){$e_2$}
\put(62.00,31.5){$e_3$}
\put(52,10.00){$e_4$}
\put(74.75,10.00){$e_5$}
\put(131.25,40){$e_1$}
\put(156.75,39.5){$e_2$}
\put(145,30){$e_3$}
\put(130.5,18.5){$e_4$}
\put(157.25,17.75){$e_5$}
\put(143,2){$C(G)$}
\end{picture}
 %%

%%%%%%%%end of input
  \caption{Line graph $L(G)$ and 
clique-inserted graph $C(G)$}
\relabel{f13}
\end{figure}

From observation (iii) above, 
$C(G)$ is in a type of connected graphs
with a matching whose removal yields 
components which are all complete graphs.
As $t(L(G))=|\setst_{C(G)}(M)|$ holds for any 
connected graph $G$ with $M$ defined above, 
we now extend our problem to
finding an expression for $|\setst_Q(M)|$, where 
$Q$ is an arbitrary connected graph and 
$M$ is any matching of $Q$
such that all components of 
$Q-M$ are complete graphs.

Throughout this section, we 
assume 
\begin{enumerate}
\item 
$Q$ is a simple and connected graph 
with a matching $M$
such that 
all components $Q_1, Q_2, \cdots, Q_n$
of $Q-M$ are complete graphs; 
% and each edge of $M$ joins two vertices in distinct components of $G-M$;
\item for $i=1,2,\cdots,n$, 
$V_i=V(Q_i)=\{v_{i,j}: j=1,2,\cdots,k_i\}$, 
where $k_i=|V_i|$; 
\item $M=\{e_1,e_2,\cdots, e_m\}$ and 
$M_i$ is the set of those edges of 
$M$ which have one end in $V_i$
and $m_i=|M_i|$ for $i=1,2,\cdots,n$;
\item 
$v_{i,j}$ is incident with an edge of $M_i$
if and only if $1\le j\le m_i$;
%Thus, 
\item 
$Q^*$ is the graph obtained from $Q$ by 
contracting all edges of $Q_i$ for all $i=1,2,\cdots,n$.
Thus each $Q_i$ is converted to a vertex in $Q^*$ denoted by $v_i$.

%Thus 
\end{enumerate}

With the above assumptions, 
we observe that
 $V(Q^*)=\{v_1, v_2,\cdots,v_n\}$ and $E(Q^*)=M$.
As $M$ is a matching of $Q$ and $Q$ is connected, we have $1\le m_i\le k_i$.
If $k_i>m_i$, then vertex $v_{i,j}$ is not incident with any edge of $M$ for all $j: m_i<j\le k_i$.
If $k_i=m_i$ for all 
$i=1,2,\cdots,n$, 
then  $|\setst_Q(M)|=t(L(Q^*))$.
Thus result (\ref{eq-r0}) is a special case 
of Theorem~\ref{gen-result} which is the main result 
to be established in this section. 

\begin{theo}\relabel{gen-result}
For $Q, Q^*$ and $M$ defined above, we have 
\begin{equation}\relabel{eq3-04}
|\setst_Q(M)|=\sum_{T\in \sett(Q^*)}
\sum_{f\in \Gamma(E(Q^*)-E(T))}
\prod_{i=1}^n k_i^{k_i-2-|f^{-1}(v_i)|}.
\end{equation}
\end{theo}

To prove Theorem~\ref{gen-result},
by the following result, 
we only need to consider the case that $k_i>m_i$ 
for all $i=1,2,\cdots,n$.

\begin{pro}\relabel{case-needed}
Theorem~\ref{gen-result} holds if it %equality (\ref{eq3-04})  
holds whenever $k_i>m_i$ for all $i=1,2,\cdots,n$.
\end{pro}

\proof Assume that $M$ is fixed and so all $m_i$'s are fixed. 
Without loss of generality, we only need to show that
with $k_i$, where $k_i\ge m_i$, to be fixed 
for all $i=2,\cdots,n$, 
if (\ref{eq3-04}) holds for every integer $k_1$ 
with $k_1\ge m_1+1$, 
then it also holds for the case $k_1=m_1$.

%Now also assume that $k_i$'s are fixed for all $i=2,3,\cdots,n$.
For any integer $k_1\ge m_1$, let 
$$
\gamma(k_1)=|\setst_Q(M)|.
$$
By the assumption, for any $k_1\ge m_1+1$,
(\ref{eq3-04}) holds and thus
\begin{equation}\relabel{eq3-05}
\gamma(k_1)=\sum_{T\in \sett(Q^*)}
\sum_{f\in \Gamma(E(Q^*)-E(T))}
k_1^{k_1-2-|f^{-1}(v_1)|}
\prod_{i=2}^n k_i^{k_i-2-|f^{-1}(v_i)|}
=\sum_{s=0}^{m_1-1} a_s k_1^{k_1-2-s},
\end{equation}
where 
\begin{equation}\relabel{eq3-05-1}
a_s=\sum_{T\in \sett(Q^*)}
\sum_{f\in \Gamma(E(Q^*)-E(T))\atop |f^{-1}(v_1)|=s}
\prod_{i=2}^n k_i^{k_i-2-|f^{-1}(v_i)|}.
\end{equation}
It is clear that $a_s$ is independent of the value of $k_1$.

Now let $Q'$ be the graph 
$Q-E(Q_1)-\{v_{1,j}: m_1<j\le k_1\}$.
So $Q'$ is independent of $k_1$.
Note that for every $T\in \setst_Q(M)$, 
$F=T-E(T[V_1])-\{v_{1,j}: m_1<j\le k_1\}$ 
is a member of $\setsf_{Q'}(M)$,
i.e., a spanning forest of $Q'$ containing all edges of $M$,
since $v_{1,j}$ is not incident with any edge of $M$
 for all $j: m_1<j\le k_1$.
Thus $\setst_Q(M)$ can be partitioned into 
$$
\setst_Q(M)=\bigcup_{F\in \setsf_{Q'}(M)}\setst_{Q''}(F),
$$
where $Q''=Q[E(F)\cup E(Q_1)]$.
It is possible that $\setst_{Q''}(F)=\emptyset$
for some $F\in \setsf_{Q'}(M)$.
But $\setst_{Q''}(F')\cap \setst_{Q''}(F'')=\emptyset$ 
for distinct $F',F''\in \setsf_{Q'}(M)$,
implying that for any $k_1=|V_1|\ge m_1$,  
$$
\gamma(k_1)=\sum_{F\in \setsf_{Q'}(M)}|\setst_{Q''}(F)|.
$$
%where $Q''$ is the spanning subgraph of $Q$ with edge set$E(F)\cup E(Q_1)$.
By Proposition~\ref{forest-stree},
for any $F\in \setsf_{Q'}(M)$,
if $F/\{v_{1,j}: 1\le j\le m_1\}$ is connected, 
then 
$$
|\setst_{Q''}(F)|
=k_1^{k_1-2+c(F-V_1)-m_1}\prod_{j=1}^{c(F-V_1)} 
|E_{F}(V_1, V(F_j))|,
$$
where $F_1, F_2, \cdots, F_{c(F-V_1)}$ are the components of 
$F-V_1$. % and $t=c(F-V_1)$.
Let $\setsf^c_{Q'}(M)$ denote the set of those 
$F\in \setsf_{Q'}(M)$ such that 
$F/\{v_{1,j}: 1\le j\le m_1\}$ is connected.
Thus, for any $k_1\ge m_1$, we have 
\begin{eqnarray}\relabel{eq3-06}
\gamma(k_1)
&=&\sum_{F\in \setsf^c_{Q'}(M)}
k_1^{k_1-2+c(F-V_1)-m_1}
\prod_{j=1}^{c(F-V_1)} |E_{F}(V_1, V(F_j))|
\nonumber
\\
&=&\sum_{s=0}^{m_1-1}b_sk_1^{k_1-2-s},
\end{eqnarray}
where 
\begin{equation}\relabel{eq3-06-1}
b_s=\sum_{F\in \setsf^c_{Q'}(M)\atop c(F-V_1)=m_1-s}
\prod_{j=1}^{c(F-V_1)} |E_{F}(V_1, V(F_j))|.
\end{equation}
As $Q'$ is independent of $k_1$, 
for any $F\in \setsf^c_{Q'}(M)$,
the expression $\prod_{j=1}^{c(F-V_1)} |E_{F}(V_1, V(F_j))|$ 
is independent of $k_1=|V_1|$
and hence $b_s$ is independent of $k_1$.

By (\ref{eq3-05}) and (\ref{eq3-06}), 
for every integer $k_1$ with $k_1\ge m_1+1$, 
we have 
\begin{equation}\relabel{eq3-07}
\sum_{s=0}^{m_1-1}a_sk_1^{k_1-2-s}=
\sum_{s=0}^{m_1-1}b_sk_1^{k_1-2-s},
\end{equation}
where $a_s$ and $b_s$ are independent of $k_1$
for all $s=0,1,2,\cdots,m_1-1$.
Considering sufficiently large values of $k_1$ in (\ref{eq3-07}), 
we come to the conclusion that $a_s=b_s$ for all $s=0,1,\cdots,m_1$,
implying that 
\begin{eqnarray*}
\gamma(m_1)
&=&\sum_{s=0}^{m_1-1}b_sm_1^{m_1-2-s}
=\sum_{s=0}^{m_1-1}a_sm_1^{m_1-2-s}\\
&=&
\sum_{T\in \sett(Q^*)}
\sum_{f\in \Gamma(E(Q^*)-E(T))}
m_1^{m_1-2-|f^{-1}(v_1)|}
\prod_{i=2}^n k_i^{k_i-2-|f^{-1}(v_i)|},
\end{eqnarray*}
implying that (\ref{eq3-04}) holds for $k_1=m_1$.
Hence the result holds. 
\proofend

In the remainder of this section, we assume that 
$k_i\ge m_i+1$ for all $i$ with $1\le i\le n$.
Thus vertex $v_{i,k_i}$ is not incident with any edge of $M$
for each $i$.
We will complete the proof of Theorem~\ref{gen-result}
by the approach explained in the two steps below: 

%then can be proved by the two steps below: 

(a) $\setst_Q(M)$ will be partitioned into 
$t(Q^*)2^{m-n+1}$ subsets denoted by $\Delta(T_0,f)$'s,
corresponding to $t(Q^*)2^{m-n+1}$ 
ordered pairs $(T_0, f)$, where $T_0\in \sett(Q^*)$
and $f\in \Gamma(E(Q^*)-E(T_0))$;

(b) then we show that for any given 
$T_0\in \sett(Q^*)$ and
$f\in \Gamma(E(Q^*)-E(T_0))$,
$$
|\Delta (T_0,f)|=\prod_{i=1}^n k_i^{k_i-2-|f^{-1}(v_i)|}.
$$

Step (a) above will be done by Algorithm B below
which determines   a spanning tree 
$T_0$ of $Q^*$ and 
a mapping $f\in \Gamma(E(Q^*)-E(T_0))$
for any given $T\in \setst_Q(M)$.

{\bf Algorithm B} ($T \in \setst_Q(M)$):

\begin{enumerate}
\item[Step B1.] Let $T_n$ be $T$;

\item[Step B2.] for $i=n,n-1,\cdots,1$, 
let $D_i=E_{T_i}(V_i)-\Phi(T_i,V_i,v_{i,k_i})$
and $T_{i-1}$ be the graph obtained from $T_i$ 
by deleting all edges in  $D_i\cup E(T_i[V_i])$ 
and identifying all vertices of $V_i$ as one, 
denoted by $v_i$, which is a vertex of $Q^*$;

\item[Step B3.] output $T_0$ and $f$,
where $f$ is a mapping from $ D_1\cup D_2\cup \cdots \cup D_n$ to $V(Q^*)$ 
defined by $f(e)=v_i$ whenever $e\in D_i$.
\end{enumerate}

By Lemma~\ref{keep-t-edges}, each graph $T_i$ produced  
in the process of running Algorithm B is indeed a tree
and thus $T_0$ is a tree in $\sett(Q^*)$.
It is also clear that 
$D_1\cup D_2\cup \cdots \cup D_n=E(Q^*)-E(T_0)$
and so the mapping $f$ output by Algorithm B 
belongs to $\Gamma(E(Q^*)-E(T_0))$.

An example is presented below. 
Let $T$ be a tree in $\setst_Q(M)$ as shown in Figure~\ref{f6}(a),
where $Q$ is a connected graph with a matching 
$M=\{e_1,e_2,\cdots, e_8\}$  
such that $Q-M$ has four components $Q_1, Q_2, Q_3$ and $Q_4$
isomorphic to complete graphs 
of orders $5,4,6,5$ respectively. 
If we run Algorithm B with this tree $T$ as its input,
then we have $T_3, T_2, T_1$ and $T_0$ as shown in 
Figure~\ref{f6} and thus 
$$
D_4=\{e_4\},
D_3=\{e_1,e_2\},
D_2=\{e_5,e_7\},
D_1=\emptyset,
$$
implying that the mapping $f\in \Gamma(E(Q^*)-E(T_0))$
output by Algorithm B,
where $E(Q^*)-E(T_0)=\{e_1,e_2,e_4,e_5,e_7\}$,
is the one given below: 
$$
f(e_1)=f(e_2)=v_3, f(e_4)=v_4,f(e_5)=f(e_7)=v_2.
$$

\begin{figure}[htb]
  \centering
\unitlength .7mm % = 1.992pt
\linethickness{0.4pt}
\ifx\plotpoint\undefined\newsavebox{\plotpoint}\fi % GNUPLOT compatibility

%%
%%%%%%%%%%%End of input
  \caption{$T\in \setst_Q(M)$ (i.e., $T_4$)
  and $T_3,T_2,T_1,T_0$}
\relabel{f6}
\end{figure}

Let $\psi$ be a mapping from $\setst_Q(M)$ to the following set
of ordered pair $(T_0,f)$'s:
$$
\{(T_0,f): T_0\in \sett(Q^*), 
f\in \Gamma(E(Q)-E(T_0))\},
$$
defined by $\psi(T)=(T_0,f)$ if $T_0$ and $f$ are 
output by running Algorithm $B$ with input $T$.
For any $T_0\in \sett(Q^*)$ and $f\in \Gamma(E(Q)-E(T_0))$,
let $\Delta(T_0,f)=\psi^{-1}(T_0,f)$.
Thus $\setst_Q(M)$  is partitioned 
into $t(Q^*)2^{m-n+1}$ subsets $\Delta(T_0,f)$'s, 
where $T_0\in \sett(Q^*)$ and $f\in \Gamma(E(Q)-E(T_0))$.

The proof of Theorem~\ref{gen-result} 
%This section 
now remains to determine the size of $\Delta(T_0,f)$ below.

\begin{pro}\relabel{no-trees}
For any $T_0\in \sett(Q^*)$ and $f\in \Gamma(E(Q^*)-E(T_0))$,
we have 
$$
|\Delta(T_0,f)|=\prod_{i=1}^n k_i^{k_i-2-|f^{-1}(v_i)|}.
$$
\end{pro}

\proof Let $D_i=f^{-1}(v_i)=\{e\in M-E(T_0): f(e)=v_i\}$ for $i=1,2,\cdots,n$.
So $D_i\subseteq M_i$.
%each edge of $D_i$ is incident with some vertex in $V_i$.
By Algorithm $B$, $T$ is a member of $\Delta(T_0,f)$
if and only if there exist trees $T_1, T_2, \cdots, T_{n-1}$ 
such that %$T_n$ is the tree $T$ and 
for $i=n,n-1,\cdots, 1$, 
the following properties hold, where $T_n$ is the tree $T$:
\begin{enumerate}
\item [(P1)] $V(T_i)=(V(T_{i-1})-\{v_i\})\cup V_i$;
%(T_{i-1})=(V(T_i)-V_i)\cup \{v_i\}$;

\item [(P2)] $T_i-V_i$ and $T_{i-1}-v_i$ are the same graph;
and
 
\item [(P3)] 
$E_{T_{i-1}}(v_i)=\Phi(T_i, V_i, v_{i,k_i})
=E_{T_i}(V_i,V(T_i)-V_i)-D_i$
and $D_i\subseteq E_{T_i}(V_i,V(T_i)-V_i)$.
\end{enumerate}

Let $U_i=\bigcup\limits_{1\le j\le i}V_j\cup
\{v_{i+1},\cdots, v_n\}$.
Observe that if properties (P1), (P2) and (P3) hold
for all $i$ with $1\le i\le n$, % for $T_i$ and $T_{i-1}$,
then $V(T_{i})=U_{i}$ for all $i=0,1,\cdots,n$.

Now let $\Delta_0=\{T_0\}$. 
Define sets $\Delta_1, \Delta_2,\cdots,\Delta_n$
as follows.
For $i=1,2,\cdots,n$, let 
$$
\Delta_i=\bigcup_{T_{i-1}\in \Delta_{i-1}} \Psi(T_{i-1}),
$$
where $\Psi(T_{i-1})$ is 
the set of all those spanning trees $T_i$ 
of $H_i$ such that 
properties (P1), (P2) and (P3) hold for $T_i$ and $T_{i-1}$
and $H_i$ is the graph with $V(H_i)=U_i$
such that $V_i$ is a clique of $H_i$,
$H_i-V_i$ is the same as $T_{i-1}-v_i$
and 
$E_{H_i}(V_i)=E_{T_{i-1}}(v_i)\cup D_i$.
Note that for each edge $e\in E_{H_i}(V_i)$,
$e$ is actually also an edge in $Q$ and we 
assume that 
$e$ joins the same pair of vertices as it does in $Q$ 
unless $e$ as an edge of $Q$ has one end in some $V_j$ with $j>i$,
while in this case this end of $e$ in $H_i$ is $v_j$.

By (P1), (P2) and (P3), $T_{i-1}$ is 
uniquely determined by any $T_i\in \Psi(T_{i-1})$.
Thus $\Psi(T'_{i-1})\cap \Psi(T''_{i-1})=\emptyset$
for any distinct members $T'_{i-1}$ and $T''_{i-1}$
of $\Delta_{i-1}$.
For any $T_{i-1}\in \Delta_{i-1}$,
observe that $\Psi(T_{i-1})$ is actually the set 
$\setst_{H_i}(V_i,T_{i-1},D_i,v_{i,k_i})$, and thus 
by Proposition~\ref{no-trees-original}, %for any $T_{i-1}\in \Delta_{i-1}$,
we have 
$$
|\Psi(T_{i-1})|=k_i^{k_i-2-|D_i|}.
$$
Hence $|\Delta_i|=k_i^{k_i-2-|D_i|}|\Delta_{i-1}|$
for all $i=1,2,\cdots,n$.
As  $\Delta(T_0,f)=\Delta_n$, the result holds.
\proofend

We end this section with a proof of Theorem~\ref{gen-result}.
 
\vspace{0.3 cm}

\noindent {\it Proof of Theorem~\ref{gen-result}}:
By Proposition~\ref{case-needed}, we may assume that 
$k_i>m_i$ for all $i=1,2,\cdots,n$.
By %Proposition~\ref{run-alg} and 
the definition of $\psi$ and 
$\Delta(T_0,f)=\psi^{-1}(T_0,f)$,
%$\Delta(T_0,f)$ and $\setst_Q(M)$, 
we have 
$$
\setst_Q(M)=\bigcup_{T_0\in \sett(H)\atop f\in \Delta(E(H)-E(T_0))} \Delta(T_0,f),
$$
where 
%As $\Delta(T_0,f)=\psi^{-1}(T_0,f)$,
the  union gives a partition of $\setst_Q(M)$. 
Thus 
$$
|\setst_Q(M)|=
\sum_{T_0\in \sett(H)\atop f\in \Gamma(E(H)-E(T_0))} 
|\Delta(T_0,f)|
=\sum_{T_0\in \sett(H)\atop f\in \Gamma(E(H)-E(T_0))} 
\prod_{i=1}^n k_i^{k_i-2-|f^{-1}(v_i)|},
$$
where the last step follows from Proposition~\ref{no-trees}.
Hence Theorem~\ref{gen-result} holds.
\proofend

\resection{Proving Theorem~\ref{gen-Sr} for $r\ge 1$}

In this section, we shall prove 
Theorem~\ref{gen-Sr} for the case $r\ge 1$.

For any graph $G$ and edge $e$ in $G$,
let $G-e$ and $G/e$ be the graphs obtained from $G$ by
deleting $e$ and contracting $e$ respectively.
The following result is  obvious.

\begin{lem}[\cite{Biggs93,BM76}]\relabel{le0}
For any graph $G$ and edge $e$ in $G$, we have
$$
t(G)=t(G-e)+t(G/e).
$$
In particular, if $e$ is a bridge of $G$,
then $t(G)=t(G/e)$.
\end{lem}

For any edge $e$ in $G$, let $G_{\bullet e}$ be the
graph obtained from $G$ by inserting a vertex on $e$
and $G_{-e}$ be the graph obtained from $G-e$
by attaching a pendent edge to each end of $e$,
as shown in Figure~\ref{f10}.
Similarly, for any $E'\subseteq E(G)$,
let $G_{\bullet E'}$ be the
graph obtained from $G$ by inserting a vertex on each edge of $E'$
and $G_{-E'}$ be the graph obtained from $G-E'$ by
attaching a pendent edge to each end of $e$ for all $e\in E'$.
Clearly $G_{\bullet E'}$  is the graph $S(G)$ when $E'=E(G)$.
%%%%%%%%%%%%%%%%%%%%%%%%%%%%%%%%%%%%%%%%%%
%%%%%%%%%%%%% Figure 1
%%%%%%%%%%%%%%%%%%%%%%%%%%%%%%%%%%%%%%%%%%
\begin{figure}[htbp]
  \centering
  %\scalebox{0.8}{\includegraphics{Figure1.eps}}
  %%%%%%%%%%\input f10.tex
  %TeXCAD (http://texcad.sf.net/) Picture. File: [f10.pic]. Options on following lines.
%\grade{\on}
%\emlines{\off}
%\epic{\off}
%\beziermacro{\on}
%\reduce{\on}
%\snapping{\off}
%\pvinsert{% Your \input, \def, etc. here}
%\quality{8.000}
%\graddiff{0.005}
%\snapasp{1}
%\zoom{4.0000}
\unitlength .8mm % = 2.276pt
\linethickness{0.4pt}
\ifx\plotpoint\undefined\newsavebox{\plotpoint}\fi % GNUPLOT compatibility
\begin{picture}(120,32.75)(0,0)
\put(17,21.375){\oval(26.5,22.75)[]}
\put(59.75,21.375){\oval(26.5,22.75)[]}
\put(106.75,21.125){\oval(26.5,22.75)[]}
\put(6.25,22.5){\circle*{2.236}}
\put(49,22.5){\circle*{2.236}}
\put(96,22.25){\circle*{2.236}}
\put(26,22.75){\circle*{2.236}}
\put(68.75,22.75){\circle*{2.236}}
\put(115.75,22.5){\circle*{2.236}}
\put(58.25,22.75){\circle*{2.236}}
\put(95.75,27.5){\circle*{2.236}}
\put(115.5,28.25){\circle*{2.236}}
\put(6.25,22.5){\line(1,0){20.25}}
\put(49,22.5){\line(1,0){20.25}}
\put(15.75,25.5){$e$}
\put(4.75,17.75){$1$}
\put(47.5,17.75){$1$}
\put(94.5,17.5){$1$}
\put(25,17.75){$2$}
\put(67.75,17.75){$2$}
\put(114.75,17.5){$2$}
\put(14,3.5){(a)}
\put(56.75,3.5){(b)}
\put(103.75,3.25){(c)}
\put(96,22.5){\line(0,1){4.75}}
%\emline(115.5,22.75)(115.25,28.25)
\multiput(115.5,22.75)(-.041667,.916667){6}{\line(0,1){.916667}}
%\end
\end{picture}
  %%%%%%%%%% end of inpput 
  \caption{\ (a)\ $G$ with edge $e$
  \quad (b) The graph $G_{\bullet e}$\quad (c) The graph $G_{-e}$}
  \relabel{f10}
\end{figure}
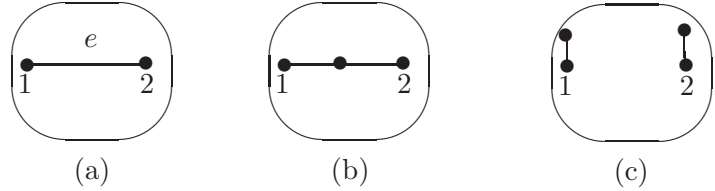

By the definition of the line graph, 
the following lemma follows from Lemma~\ref{le0}.

\begin{lem}\relabel{le1}
Let $G$ be any graph and $e$ be an edge in $G$.
Then
$$
t(L(G_{\bullet e}))
=t(L(G))+t(L(G_{- e})).
$$
In particular, if $e$ is a bridge of $G$, then
$t(L(G_{\bullet e}))=t(L(G))$.
\end{lem}

For any edge $e$ in $G$ and any non-negative integer $r$,
let $G_{r\bullet e}$ be the
graph obtained from $G$ by inserting $r$ new vertices on $e$,
i.e., replacing $e$ by a path of length $r+1$
connecting the two ends of $e$.
For any subset $F$ of $E(G)$,
let $G_{r\bullet F}$ be the
graph obtained from $G$ by replacing each edge $e$ of $F$ 
by a path of length $r+1$
connecting the two ends of $e$.

\begin{lem}\relabel{le2}
Let $G$ be any graph and $F$ be any subset of $E(G)$.
Then, for any $r\ge 0$, 
\begin{equation}\relabel{le-eq1}
t(L(G_{r\bullet F}))
=\sum_{E'\subseteq F}r^{|E'|}t(L(G_{-E'})).
\end{equation}
\end{lem}

\proof Note that for any two vertices $u,v$ in 
a graph $H$, if $N_H(u)=\{v\}$ and $d_H(v)=2$, then 
$t(L(H))=t(L(H-u))$. 
Thus, for any edge $e$ of $G$ and any positive integer $r$, 
by Lemma~\ref{le1}, we have 
\begin{equation}\relabel{le-eq2}
t(L(G_{r\bullet e}))
=t(L(G_{(r-1)\bullet e}))+t(L(G_{- e})),
\end{equation}
where $G_{0\bullet e}$ is $G$. 
Applying (\ref{le-eq2}) repeatedly deduces that 
\begin{equation}\relabel{le-eq20}
t(L(G_{r\bullet e}))
=t(L(G))+rt(L(G_{- e})).
\end{equation}
Note that (\ref{le-eq1}) is obvious for $F=\emptyset$
or $r=0$.
Now assume that $e\in F$ and $r\ge 1$. 
By induction, we have 
\begin{equation}\relabel{le-eq3}
t(L(G_{r\bullet F-\{e\}}))
=\sum_{E'\subseteq F-\{e\}}r^{|E'|}t(L(G_{-E'})).
\end{equation}
By  (\ref{le-eq20}), we have 
\begin{equation}\relabel{le-eq4}
t(L(G_{r\bullet F}))
=t(L(G_{r\bullet F-\{e\}}))+rt(L((G_{r\bullet F-\{e\}})_{-e})).
\end{equation}
Thus (\ref{le-eq1}) follows immediately from 
(\ref{le-eq3}).
\proofend

We are now ready to prove Theorem~\ref{gen-Sr} for 
the case $r\ge 1$.

\vspace {0.3 cm}

\noindent{\it Proof of Theorem~\ref{gen-Sr} for $r\ge 1$}: 
%The result for the case $r=0$ has been proved in Section 3. 
Assume that $r\ge 1$.
By Lemma~\ref{le2}, we have
\begin{equation}\relabel{theo-proof-eq1}
t(L(S_r(G)))
=\sum_{E'\subseteq E(G)}r^{|E'|}t(L(G_{-E'})).
\end{equation}
%Note that $t(L(G_{-E'}))=0$ if $E'$ contains any bridge of $G$.
The above summation needs only to take those 
subsets $E'$ of $E(G)$ with 
$t(L(G_{-E'}))>0$ (i.e. $G-E'$ is connected).
Now let $E'$ be any fixed subset of $E(G)$ 
such that $G-E'$ is connected and 
let $H$ denote $G_{-E'}$.
By Theorem~\ref{gen-Sr} for $r=0$ 
(i.e., (\ref{eq-r0})),  
\begin{equation}\relabel{theo-proof-eq2}
t(L(G_{-E'}))
=\sum_{T'\in \sett(H)}
\sum_{g\in \Gamma(E(H)-E(T'))} 
\prod_{v\in V(H)} d_{H}(v)^{d_{H}(v)-2-|g^{-1}(v)|}.
\end{equation}
Observe that $V(G)\subseteq V(H)$.
For any $v\in V(H)$,
if $v\in V(G)$, then $d_{H}(v)=d_G(v)$;
otherwise, $d_{H}(v)=1$. %it is of degree $1$ in $H$.
Thus 
\begin{equation}\relabel{theo-proof-eq3}
\prod_{v\in V(H)} d_{H}(v)^{d_{H}(v)-2-|g^{-1}(v)|}
=\prod_{v\in V(G)} d_{G}(v)^{d_{G}(v)-2-|g^{-1}(v)|}.
\end{equation}
For each $T'\in \sett(H)$, $T'$ contains all pendent edges
in $H$
and so $T'$ corresponds to $T$, where $T=T'[V(G)]$,
which is a spanning tree of $G-E'$.
 % is the subgraph of $T'$ induced by $V(G)$.
Thus $E(H)-E(T')=E(G-E')-E(T)$ and 
\begin{equation}\relabel{theo-proof-eq4}
t(L(G_{-E'}))
=\sum_{T\in \sett(G-E')}
\sum_{g\in \Gamma(E(G-E')-E(T))} 
\prod_{v\in V(G)} d_{G}(v)^{d_{G}(v)-2-|g^{-1}(v)|}.
\end{equation}
%implying that 
By (\ref{theo-proof-eq1}) and (\ref{theo-proof-eq4}),
\begin{equation}\relabel{theo-proof-eq5}
t(L(S_r(G)))=\sum_{E'\in E(G)}r^{|E'|}
\sum_{T\in \sett(G-E')}
\sum_{g\in \Gamma(E(G-E')-E(T))} 
\prod_{v\in V(G)} d_{G}(v)^{d_{G}(v)-2-|g^{-1}(v)|}.
\end{equation}
By replacing $E(G)-E'-E(T)$ by $E''$,  
(\ref{theo-proof-eq5}) implies that 
\begin{eqnarray*}
t(L(S_r(G)))
&=&\sum_{E''\subseteq E(G)}
\sum_{T'\in \sett(G-E'')}
\sum_{g\in \Gamma(E'')} 
r^{|E(G)|-|E''|-|E(T')|} 
\prod_{v\in V(G)} 
d_{G}(v)^{d_{G}(v)-2-|g^{-1}(v)|}\\
&=&\sum_{E''\subseteq E(G)}
r^{|E(G)|-|E''|-|V(G)|+1} t(G-E'')
\sum_{g\in \Gamma(E'')} 
\prod_{v\in V(G)} 
d_{G}(v)^{d_{G}(v)-2-|g^{-1}(v)|}\\
&=&\sum_{E'''\subseteq E(G)}
r^{|E'''|-|V(G)|+1} t(G[E'''])
\sum_{g\in \Gamma(E-E''')} 
\prod_{v\in V(G)} 
d_{G}(v)^{d_{G}(v)-2-|g^{-1}(v)|}.
\end{eqnarray*}
Hence the case $r\ge 1$ of Theorem~\ref{gen-Sr} holds.
\proofend

\section{Proof of Conjecture~\ref{con-yan}}

Now we turn back to those connected graphs $G$
mentioned in Conjecture~\ref{con-yan} and 
apply the following result and Theorem~\ref{gen-Sr}
to deduce a relation between $t(L(S_r(G)))$ and $t(G)$.
The case $r=1$ of this relation 
is exactly the conclusion of Conjecture~\ref{con-yan}.

\begin{lem}\relabel{equality-spanning trees}
Let $H$ be any connected graph of order $n$ and size $m$.
For any integer $i$ with $0\le i\le m-n+1$,
we have
$$
{m-n+1\choose i}t(H)=
\sum_{E'\subseteq E(H)\atop |E'|=i}t(H-E').
$$
%where $H-E'$ is the spanning subgraph of $H$ with edge set $E(H)-E'$.
\end{lem}

\proof We prove this result by 
providing two different methods to 
determining the size of the following set:
$$
\Theta=\{(T,E'): T \mbox{ is a spanning tree of }H
\mbox{ and }
E'\subseteq E(H)-E(T) \mbox{ with }|E'|=i\}.
$$
Note that for each spanning tree $T$ of $H$, 
as $|E(H)|=m$ and $|E(T)|=n-1$,
the number of subsets $E'$ of $E(H)-E(T)$ with 
$|E'|=i$ is ${m-n+1\choose i}$.
On the other hand, for each $E'\subseteq E(H)$ with
$|E'|=i$, there are exactly $t(H-E')$ spanning trees $T$ of $G$ 
such that
$E'\subseteq E(H)-E(T)$. Thus the result holds. 
\proofend

We now deduce the following consequence of 
Theorem~\ref{gen-Sr} for those connected graphs $G$
mentioned in Conjecture~\ref{con-yan}.

\begin{cor}\relabel{result-extend}
Let $G$ be a connected graph of order $n+s$ and size $m+s$ 
in which $s$ vertices are of degree $1$ and all others 
are of degree $k$, where $k\ge 2$.
Then, for any $r\ge 0$, 
$$
t(L(S_r(G)))=k^{m+s-n-1}(rk+2)^{m-n+1}t(G).
$$
\end{cor}

\proof For any $E'\subseteq E(G)$ with $t(G[E'])\ne 0$,
$E'$ contains every bridge of $G$, and so
$d(u_e)=d(v_e)=k$ for all $e\in E(G)-E'$.
By Theorem~\ref{gen-Sr}, we have 
\begin{eqnarray*}
t(L(S_r(G)))&=&(k^{k-2})^n\sum_{E'\subseteq E(G)}
t(G[E'])r^{|E'|-(n+s)+1} (2k^{-1})^{(m+s)-|E'|}\\
&=&(k^{k-2})^n
r^{-(n+s)+1} (2k^{-1})^{(m+s)}
\sum_{E'\subseteq E(G)}
t(G[E'])r^{|E'|} (2k^{-1})^{-|E'|}\\
&=&(k^{k-2})^n
r^{-(n+s)+1} (2k^{-1})^{(m+s)}
\sum_{E''\subseteq E(G)}
t(G-E'')r^{|E(G)|-|E''|} (2k^{-1})^{|E''|-|E(G)|}\\
&=&(k^{k-2})^n
r^{-(n+s)+1} (2k^{-1})^{(m+s)}
\sum_{j=0}^{m-n+1}r^{m+s-j} (2k^{-1})^{j-m-s}
\sum_{E''\subseteq E(G)\atop |E''|=j}
t(G-E'')\\
&=&(k^{k-2})^n %r^{-(n+s)+1} 
\sum_{j=0}^{m-n+1}r^{m-n+1-j} (2k^{-1})^{j}
{m-n+1\choose j} t(G)\hspace{1.3 cm} 
(\mbox{by Lemma~\ref{equality-spanning trees}}) \\
&=&(k^{k-2})^n (r+2k^{-1})^{m-n+1}t(G) \\
&=&k^{n(k-2)-(m-n+1)}(kr+2)^{m-n+1}t(G) \\
&=&k^{m+s-n-1}(kr+2)^{m-n+1}t(G),
\end{eqnarray*}
where 
%the fourth last expression follows from Lemma~\ref{equality-spanning trees}and 
the last expression follows from the equality 
$2(m+s)=kn+s$ by the given conditions on $G$.
Hence the result is obtained.
\proofend

Notice that (\ref{eq1-3}) is the special case of 
Corollary~\ref{result-extend} for $r=0$ 
while the conclusion of 
Conjecture~\ref{con-yan} is the special case of 
Corollary~\ref{result-extend} for $r=1$.

We end this section with the 
the following result on some special bipartite graphs,
which can be obtained by applying 
Lemma~\ref{equality-spanning trees}
and the case $r=0$ of Theorem~\ref{gen-Sr}.

\begin{cor}\relabel{bipar}
Let $G=(A,B;E)$ be a connected bipartite graph of order $n$ and size $m$ 
such that $d(x)\in \{1,a\}$ for all $x\in A$ and 
$d(y)\in \{1,b\}$ for all $y\in B$, where $a\ge 2$ and $b\ge 2$.
Then
$$
t(L(G))=a^{(a-2)n_1}b^{(b-2)n_2}(a^{-1}+b^{-1})^{m-n+1}t(G),
$$
where $n_1$ is the number of vertices $x$ in $A$ with  $d(x)=a$
and  $n_2$ is the number of vertices $y$ in $B$ with $d(y)=b$.
\end{cor}

%Corollary~\ref{bipar} can be proved similarly as Corollary~\ref{result-extend} by applying Lemma~\ref{equality-spanning trees}and the case $r=0$ of Theorem~\ref{gen-Sr}.

The result of Corollary~\ref{bipar} in the case that
$G$ is an $(a, b)$-semiregular bipartite graph
was originally due to 
Cvetkovi\'{c} (see Theorem 3.9 in \cite{mac}, \S 5.2 of \cite{moh}, 
or \cite {sat}).

\vspace{0.5 cm}

\noindent {\it Acknowledgement.} 
The authors wish to thank the referees 
for their very helpful suggestions.


\begin{thebibliography}{99}


\bibitem{Aig} %Aigner, Martin and Ziegler, G\"{u}nter M.,
M. Aigner and G. Ziegler,
{\it Proofs from The Book},
Fourth edition. 
Springer-Verlag, Berlin, 2010. 
%viii+274 pp. 

\bibitem{ber}
A. Berget, A. Manion, M. Maxwell, A. Potechin, V. Reiner, 
The critical group of a line graph,
{\it Ann. Comb.} {\bf 16} (2012), 449-488.

\bibitem{bid}
H. Bidkhori and S. Kishore, Counting spanning trees of a directed line graph, arXiv: 0910.3442v1.

\bibitem{Biggs93}
N. L. Biggs, \textit{Algebraic Graph Theory}, 2nd edn, Cambridge,
Cambridge University Press, 1993.

\bibitem{BM76}
J. A. Bondy and U. S. R. Murty, \textit{Graph Theory with
Applications}, American Elsevier, New York, 1976.


\bibitem{che}
H. Y. Chen, F. J. Zhang, 
The critical group of a clique-inserted graph, 
{\it Discrete Math.} {\bf 319} (2014), 
24-32.

\bibitem{cve}
D.Cvetkovi\'{c}, M.Doob, H.Sachs,
Spectra of Graphs.
Theory and Application,
Pure Appl. Math.,vol. 87,
Academic Press,Inc.
[Harcourt Brace Jovanovich, Publishers], New York, London,1980.


\bibitem{kel} A.K.Kelmans,
On properties of the characteristic polynomial of a graph,
in:Kibernetiku Na Sluzbu Kom., vol.4, Gosener-goizdat,
Moscow, 1967, pp.27-47(in Russian).

\bibitem{knu}
D.E. Knuth. Oriented subtrees of an arc digraph, 
{\it J. of Combin. Theory} {\bf 3} (1967), 309-314.

\bibitem{lev}
L. Levine, Sandpile groups and spanning trees of directed line graphs,
{\it J. of Combin. Theory Ser. A }
{\bf 118} (2011), 350-364.

\bibitem{lov}
L. Lov\'{a}sz,
{\it Combinatorial Problems and Exercises},
North-Holland, Amsterdam (1979).

\bibitem{lov2} L. Lov\'{a}sz, M.D. Plummer, 
Matching Theory, 
{\it Ann. Discrete Math.} {\bf 29}, 
North-Holland, Amsterdam, 1986.


\bibitem{mac}
I.G. Macdonald, 
{\it Symmetric functions and Hall polynomials}, 
2nd edition. Oxford Mathematical Monographs. 
Oxford Science Publications. 
The Clarendon Press, Oxford University Press, New
York, 1995.


\bibitem{mad}
W. Mader, 
Minimale $n$-fach kantenzusammenhangende Graphen, 
{\it Math. Ann.} {\bf 191} (1971), 21-28.

%\bibitem{Mey2011}Aaron Meyerowitz,  http://mathoverflow.net/questions/66588/number-of-spanning-forests-in-a-graph.

\bibitem{moh} B. Mohar, 
The Laplacian Spectrum of Graphs, 
{\it Graph Theory, Combinatorics, and Applications} {\bf 2}
Ed. by Y. Alavi, G. Chartrand, 
O. R. Oellermann, A. J. Schwenk. Wiley, 1991, 871-898. 



\bibitem{per}
D. Perkinson, N. Salter, T. Y. Xu, 
A note on the critical group of a line graph, 
{\it Electron. J. Combin.} {\bf 18} (2011), \#P124.

\bibitem{sat}
I. Sato, 
Zeta functions and complexities of a semiregular bipartite graph and 
its line graph, 
{\it Discrete Math.} {\bf 307} (2007), 237-245. 


\bibitem{shi}
T. Shirai, 
The spectrum of infinite regular line graphs, 
{\it Trans. Amer. Math. Soc.}
{\bf 352} (2000), no. 1, 115-132.


\bibitem{vah}
E.B.Vahovskii,
On the characteristic numbers of incidence matrices for non-singular graphs,
{\it Sibirsk. Mat. Zh.} {\bf 6} (1965), 44-49 (in Russian).


\bibitem{yan1} Weigen Yan,
On the number of spanning trees of some irregular line graphs,
{\it J. Combin. Theory Ser. A} {\bf 120} (2013), 1642-1648.

\bibitem{zha}
F.J.Zhang, Y.-C.Chen, Z.B.Chen,
Clique-inserted-graphs and spectral dynamics of clique-inserting,
{\it J. Math. Anal. Appl.} {\bf 349} (2009), 211-225.

\bibitem{zha2}
Z. H. Zhang, Some physical and chemical indices of clique-inserted lattices, Journal of Statistical    
   Mechanics: 
{\it Theory and Experiment}, 
doi:10.1088/1742-5468/2013/10/P10004.


\end{thebibliography}
\end{document}